\documentclass[12pt]{article}

\usepackage{amsmath}
\usepackage{amsfonts}
\usepackage{amssymb}
\usepackage{amscd}
\usepackage{amsthm}
\usepackage{latexsym}

\def\C{\mathbb{C}}
\def\FFF{\mathbb{F}}
\def\N{\mathbb{N}}
\def\R{\mathbb{R}}
\def\T{\mathbb{T}}

\def\A{{\mathcal A}}
\def\B{{\mathcal B}}
\def\BC{BC}
\def\CC{C}
\def\D{{\mathcal D}}
\def\F{\mathcal F}
\def\H{\mathcal H}
\def\K{\mathcal K}
\def\M{\mathcal M}
\def\NN{\mathcal N}
\def\S{\mathcal S}
\def\L{\mathcal L}
\def\U{\mathcal U}
\def\CCC{{\mathfrak C}}
\def\FF{{\mathfrak F}}
\def\h{{\mathfrak h}}

\def\KK{{\mathfrak K}}

\def\W{{\mathcal W}}
\def\Rep{\mathfrak{Rep}}

\def\de{\mathrm{d}}

\def\e{\varepsilon}

\def\im{\mathop{\mathsf{Im}}\nolimits}

\def\I{{\rm 1\kern-.26em I}}
\def\SA{S_{\!\A}}
\def\FA{F_{\!\A}}
\def\pF{\pi_{\!\hbox{\it \tiny F}}}
\def\iF{{\hbox{\it \tiny F}}}
\def\oB{\omega^{\hbox{\it \tiny B}}}
\def\gB{\gamma^{\hbox{\it \tiny B}}}
\def\gtB{\tilde{\gamma}^{\hbox{\it \tiny B}}}

\def\var{\mathfrak z} 
\def\Op{\mathfrak{Op}}
\def\CBA{{\mathfrak C}^{B}_{\!\A}}
\def\BBA{{\mathfrak B}^{B}_{\!\!\A}}

\def\1{\mathfrak{1}}
\def\0{\mathfrak{0}}

\newtheorem{Lemma}{Lemma}[section]
\newtheorem{Corollary}[Lemma]{Corollary}
\newtheorem{Theorem}[Lemma]{Theorem}
\newtheorem{Proposition}[Lemma]{Proposition}
\newtheorem{Definition}[Lemma]{Definition}
\newtheorem{Remark}[Lemma]{Remark}

\numberwithin{equation}{section}

\begin{document}

\title{Spectral and Propagation Results \\
for Magnetic Schr\"odinger Operators;\\
a C$^*$-Algebraic Framework}

\date{\today}				

\author{Marius M\u antoiu, Radu Purice\footnote{Institute
of Mathematics Simion Stoilow of
the Romanian Academy, P.O.  Box 1-764, Bucharest, RO-70700, Romania.
Email: mantoiu@imar.ro, purice@imar.ro}
\ and Serge Richard\footnote{Institut Girard Desargues,
Universit\'e Claude Bernard Lyon 1, 43 avenue du 11 novembre 1918,
69622 Villeurbanne cedex, France. 
Email: srichard@igd.univ-lyon1.fr}}

\maketitle
\abstract{We study generalised magnetic Schr\"odinger operators of the form
$H_{\!h\!}(A,V)=h(\Pi^A)+V$, where $h$ is an elliptic symbol, 
$\Pi^A=-i\nabla-A$, with $A$ a vector potential defining a variable 
magnetic field $B$, and $V$ is a scalar potential. 
We are mainly interested in anisotropic functions $B$ and $V$. 
The first step is to show that these operators are affiliated
to suitable $C^*$-algebras of (magnetic) pseudodifferential operators. 
A study of the quotient of these $C^*$-algebras by the ideal of compact
operators leads to formulae for the essential spectrum of $H_{\!h\!}(A,V)$,
expressed as a union of spectra of some asymptotic operators, supported
by the quasi-orbits of a suitable dynamical system. The quotient 
of the same $C^*$-algebras by other ideals give localization
results on the functional calculus of the operators $H_{\!h\!}(A,V)$, 
which can be interpreted as non-propagation properties of their unitary 
groups.}\footnote{\textbf{2000 Mathematics Subject Classification:} 35P05, 35S05, 46L55, 46L60, 47A10, 47A60, 47G30, 47L65, 81R15, 81Q10

\textbf{\ \ Keywords:} Magnetic field, pseudodifferential operators, dynamical system, twisted crossed-product, essential spectrum, Schr\"{o}dinger operators}

\section*{Introduction}

Until rather recently, the nature of the essential spectrum of self-adjoint
partial differential operators with anisotropic coefficients was poorly
understood. 
It was clear that what counts is the long-scale behaviour of these
coefficients, but it was not clear how to express this in a general and
unified manner.

In recent years a significant progress was achieved. We do not intend to trace
the history of this topic here. Let us mention however that,
in our opinion, the most efficient tools belong to the theory of $C^*$-algebras in
conjunction with some configurational framework as dynamical systems or
Lie groupoids. We refer for example to \cite{ABG}, \cite{GI1}, \cite{GI2}, 
\cite{GI3}, \cite{LN}, \cite{N}, \cite{M1}, \cite{M2} and
references therein. The list is by no means complete and we do not
describe in detail the different but connected points of view 
of these works. 
We only say some words on the common part of the ideas
involved ({\it cf}.~also \cite{GI3}).
Let $H$ be a self-adjoint operator in a Hilbert space $\H$. 
The central fact is that this operator is {\it affiliated} to some 
$C^*$-algebra $\CCC$ of bounded operators in $\H$; 
this means that its resolvent family belongs to $\CCC$. 
The essential spectrum of $H$ can be calculated if we understand the image
of $H$ in the quotient $C^*$-algebra $\CCC/\CCC\cap K(\H)$, 
where $K(\H)$ is the ideal of compact operators in $\H$. 
Many choices for $\CCC$ are possible;
the skill consists in choosing one for which the quotient is comprehensible. 
This is usually done by keeping track of the subjacent configuration space $X$
of the problem. Such a space is available because we study differential, 
or more generally, pseudodifferential operators. In
\cite{GI2}, \cite{GI3}, \cite{M1}, \cite{M2} it is assumed that $X$ is an abelian, 
locally compact group, and this leads to a dynamical system background for the 
$C^*$-algebras. In \cite{LN} and \cite{N} the authors work in a smooth groupoid
setting, which is very general.

All these mentioned works leave aside {\it variable} magnetic fields and 
this is the topic we address in the present article. 
Spectral analysis for magnetic (pseudo-)differential operators may be 
considered as a difficult matter. 
One of the reasons is gauge covariance: The vector potential $A$ defining the 
magnetic field $B$ by the relation $B=dA$ and appearing in the explicit 
expression of the operator is highly non-unique and largely irrelevant.
What counts is the magnetic field, which is hidden. 
Obviously, good spectral results should be expressed with no reference to any
vector potential. 
On the other hand, it is clear that the magnetic field plays a very 
different role than a scalar potential. 
Thus, one needs $C^*$-algebras incorporating naturally
magnetic fields, in a manifestly invariant way.

Based on works done in \cite{MP1}, \cite{MP2}, it was shown in $\cite{MPR}$
how to achieve this. 
The key concept is that of {\it twisted crossed product $C^*$-algebras}.
These algebras have been developed in a much more general setting in
\cite{BS}, \cite{PR1} and \cite{PR2}.
They are a more sophisticated version of the wide-known notion of crossed 
product algebras, already used in connection with anisotropic operators 
(without magnetic fields) in \cite{GI2}, \cite{GI3}, \cite{AMP} or \cite{M1}. 
We have shown in \cite{MPR} that certain twisted crossed
products are related to a twisted version of the Weyl pseudodifferential
calculus, introduced in \cite{KO1}, \cite{MP2} and \cite{KO2}, which is the
natural pseudodifferential calculus when twisted observables as magnetic
momenta are present. This is also basic for a strict deformation quantization
\`a la Rieffel for physical systems placed in magnetic fields, as explained in
\cite{MP3}.

It will be shown in the present article that twisted crossed product algebras 
and their natural Hilbert space representations are the right structures that
lead to results on the essential spectrum of magnetic Schr\"odinger
operators.
To describe briefly the output, let us consider in $\R^N$ an elliptic 
symbol $h$, a magnetic field $B$, a vector potential $A$ for the magnetic field and a scalar potential $V$. 
Let $H_{\!h\!}(A,V)$ denote the operator $h(\Pi^A)+V$.
We shall prove that the essential spectrum 
$\sigma_{\hbox{\rm \tiny ess}}\big(H_{\!h\!}(A,V)\big)$ of $H_{\!h\!}(A,V)$ is 
equal to $\overline{\cup_\nu\sigma[H_{\!h\!}(A_\nu,V_\nu)]}$, where $A_\nu$ 
is a vector potential for the magnetic field $B_\nu$. Here $B_\nu$ and $V_\nu$ 
are defined respectively by the asymptotic behaviour of the magnetic 
field $B$ and of the scalar potential $V$ at infinity. Actually this behaviour
is codified by a $C^*$-algebra of functions on $\R^N$. The Gelfand spectrum of 
this $C^*$-algebra is a compact dynamical system and the functions $B_\nu$ and
$V_\nu$ are just restrictions of $B$ and $V$ to quasi-orbits of
this dynamical system situated at infinity.

In a sense that will be discussed in Subsection \ref{VOX}, this 
result is a consistent extension, in the bounded case, of similar results of 
Helffer and Mohamed obtained by strictly analytical methods for a restricted 
class of perturbed magnetic Laplacians. 
However, in \cite{HM} the potential $V$ and the magnetic field
$B$ are allowed to be unbounded (under suitable restrictions). Unbounded 
scalar potentials $V$ are also considered in a great generality in 
\cite[Sec.~3]{GI3}, but we are mainly interested in the case when a 
magnetic field is also present. 
We stress that our results are valid for any elliptic symbol $h$, and not only 
for the usual magnetic Laplacian. 
There exist very few spectral results in the literature in such a general 
framework.

Actually, a $C^*$-algebraic setting can support other problems in the spectral
theory of self-adjoint operators than just calculating essential spectra. In
\cite{ABG} and \cite{R}, $C^*$-algebras are used in order to
get a Mourre estimate, which is basic for obtaining useful resolvent
estimates, finer spectral properties and scattering theory. 
Such developments require usually more detailed informations about the models 
under study and cannot be done for magnetic operators in the very general 
setting in which we will be placed below. 
But there is still a spectral topic that is available in the present 
generality, that of {\it localization properties}. 
Such results say roughly that if the support of a continuous function $\eta$
does not intersect the spectrum of an
asymptotic operator $H_{\!h\!}(A_\nu,V_\nu)$, then the operator
$\eta[H_{\!h\!}(A,V)]$ will be small when localized in the
neighbourhood of the quasi-orbit that defines $H_{\!h\!}(A_\nu,V_\nu)$. This
has as an immediate consequence a non-propagation statement for the unitary
group generated by $H_{\!h\!}(A,V)$: If a state has a spectral support with 
respect to $H_{\!h\!}(A,V)$ which do not intersect the spectrum of the
asymptotic operator $H_{\!h\!}(A_\nu,V_\nu)$, then this state cannot evolve 
under the unitary evolution generated by $H_{\!h\!}(A,V)$ towards the 
corresponding quasi-orbit.
We refer to Subsection \ref{subsecprop} for a precise statement, to Section 
\ref{secprop} for the proof and to \cite{AMP} and \cite{M2} for more explanations 
in the case $B=0$. Some particular examples can also be found in \cite{DS}.

Let us finally describe the content of this article.
In Section \ref{secdescription} we introduce the framework, recall some 
useful formulae and state precisely all the results mentioned above. 
A powerful affiliation criterion is exposed in Theorem \ref{thmaff} and 
Corollary \ref{corol}, the essential spectrum is calculated in Theorem 
\ref{thmess}, and propagation results are contained in Theorem \ref{thmprop}. 
Section \ref{secaffiliation} is mainly devoted to the proof of the 
affiliation criterion. It is the most technical part of this paper, 
but is also of central importance for the other results.
The ingredients used for the description of the essential spectrum are 
explained in Section \ref{secess}, and an abstract version of Theorem 
\ref{thmess} is presented and proved. 
The short proof of the propagation property is given in Section \ref{secprop}. 
And the last section is dedicated to examples and to a comparison with the
results of \cite{HM}.  

By some extra technical effort we could have obtained certain minor
ameliorations of the results. Sometimes this will be rather evident to the
attentive reader. For sake of simplicity we stick to the present version. The
main goal of an improved subsequent work would be to allow continuous,
unbounded functions $B_{jk}$ and $V$. At least for the moment we do not know
how to do that.    

\paragraph{Notations:} We briefly set together some conventions and notations.
$X$ denotes the euclidean space $\R^N$, with $N \in \N$, and $X^\star$ 
denotes its dual space, commonly identified with $\R^N$.  
The Lebesgue measures on $X$ and $X^\star$ are normalized in such a way 
that the Fourier transform $\F: L^1(X) \to \CC_0(X^\star)$, with 
$(\F f)(p)=\int_X \de x\,e^{ip \cdot x}f(x)$, induces a unitary 
map from $L^2(X)$ to $L^2(X^\star)$.
$\BC(X)$, $\BC_u(X)$ and $\CC_0(X)$ are respectively the algebra of bounded and
continuous functions on $X$, the algebra of bounded and uniformly continuous 
functions on $X$ and the ideal of continuous functions on $X$ that converge 
to $0$ at infinity. We denote by $BC^\infty(X)$ the space of smooth complex 
functions on $X$ with bounded derivatives of any order.
Except in some specified and well defined context, $\H$ is the Hilbert space
$L^2(X)$, $\B(\H)$ denotes the algebra of bounded operators in $\H$,
and $K(\H)$ the ideal of compact operators in $\H$.

\section{Main results}\label{secdescription}

\subsection{The framework}\label{fram}

In recent papers a pseudodifferential calculus (\cite{KO1}, \cite{MP2} and
\cite{KO2}) and an algebraic framework (\cite{MP1} and \cite{MPR}) where 
introduced in order to deal with the quantization problem for a particle 
in a magnetic field. 
We start by recalling very briefly some aspects of each construction. 
These approaches are complementary and both are relevant for generalised
magnetic Schr\"odinger operators and, specifically, for the statement 
of our main results. 
The relations between these approaches, rigorously investigated in 
\cite[Sec.~3.1 \& 3.2]{MPR}, will be outlined at the end of the section.
We refer to the publications quoted above for more infor\-mations and details.

\subsubsection*{The magnetic Weyl calculus}

We recall the magnetic version of the usual Weyl calculus and 
the associated magnetic symbolic calculus. 
The corresponding magnetic Moyal algebra is also introduced with a brief 
review of some of its properties.

Assume that $B$ is a continuous magnetic field on $X$, and let $A$ be a 
continuous vector potential that generates the magnetic field, 
{\it i.e.}~$A$ is a 1-form on $X$ that satisfies 
$\partial_j A_k - \partial_k A_j = B_{jk}$.
In the Hilbert space $\H$, $Q_j$ denotes the operator of multiplication by 
the $j$th coordinate, and $\Pi^A_j:=-i\partial_j -A_j$ is the $j$th component 
of the usual magnetic momentum. 
The magnetic Weyl calculus is a gauge covariant prescription that assigns to 
suitable symbols $f:X \times X^\star \to \C$ an operator $\Op^A(f)\equiv 
f(Q,\Pi^A)$ acting in $\H$.
More precisely, if $\rho$ is a scalar function on $X$ and $A':=A+\nabla \rho$
is another vector potential that generates the same magnetic field, then the 
relation $e^{i\rho}\;\!\Op^A(f)\;\!e^{-i\rho}=\Op^{A'}(f)$ holds.
The prescription is formally given, for any $u \in \H$, by  
 \begin{equation}\label{Op}
[\Op^A(f)\;\!u](x):=\int_X \de y \int_{X^\star} \de p \;e^{ip\cdot(x-y)}
\;\!\lambda^A(x;y-x)\;\!f\big(\hbox{$\frac{1}{2}$}(x+y),p\big)\;\!u(y),
\end{equation}
where
\begin{equation}\label{defdelambda}
\lambda^A(q;x):=\exp\big(-i\Gamma^A[q,q+x]\big) 
\end{equation}
and $\Gamma^A[q,q+x]$ is the circulation of $A$ along the segment of ends 
$q$ and $q+x$. 

For brevity let us denote by $\Xi$ the phase space $X \times X^\star$.
The magnetic symbolic calculus is a noncommutative composition law 
$\circ$ acting on functions $f,g: \Xi \to \C$ such that the 
relation $\Op^A(f\circ g)=\Op^A(f)\;\Op^A(g)$ is satisfied. This operation, 
called {\it the magnetic Moyal product}, is formally defined, for $\xi=(q,p)$, 
$\eta=(x,k)$ and $\zeta=(y,l)$ in $\Xi$, by 
\begin{equation}\label{produitMoyal}
[f\circ g](\xi):=4^N\int_\Xi \de\eta \int_\Xi 
\de \zeta \;\!e^{-2i\sigma(\eta,\zeta)}\;\!\oB\big(q-x-y;2x,2(y-x)\big)
\;\!f(\xi-\eta)\;\!g(\xi-\zeta),
\end{equation}
where
\begin{equation}\label{defdeomega}
\oB(q;x,y):=\exp\big(-i\Gamma^B\langle q,q+x,q+x+y\rangle\big)
\end{equation}
and $\Gamma^B\langle q,q+x,q+x+y\rangle$ is the flux of the magnetic field 
through the triangle defined by the points $q$, $q+x$ and $q+x+y$.
An explicit parametrized formula for $\oB(q;x,y)$ is given in equation
\eqref{defdeux}.  
The expression $\sigma(\eta,\zeta)$ in \eqref{produitMoyal} is equal to
$k \cdot y - l \cdot x$.
Let us mention that an involution can also be defined by 
$f^\circ(\xi):=\overline{f(\xi)}$ and satisfies $\Op^A(f^\circ)=\Op^A(f)^*$.

The integrals defining $f\circ g$ are absolutely convergent only for 
restricted classes of symbols. In order to deal with more general 
distributions, an extension by duality was proposed in \cite{MP2} under an 
additional smoothness condition on the magnetic field. So let us assume
that the components of the magnetic field are 
$\CC^\infty_{\hbox{\tiny \rm pol}}(X)$-functions, {\it i.e.}~they are 
indefinitely derivable and each derivative is polynomially bounded.
The duality approach is based on the observation \cite[Lem.~14]{MP2}\;: 
For any $f,g$ in the Schwartz space $\S(\Xi)$, we have  
\begin{equation*}
\int_\Xi \de \xi \;\![f\circ g](\xi)=\int_\Xi \de \xi \;\![g\circ f](\xi)=
\int_\Xi \de \xi \;\!f(\xi)\;\!g(\xi) 
=\langle \overline{f},g\rangle \equiv (f,g).  
\end{equation*}
As a consequence, if $f,g$ and $h$ belong to $\S(\Xi)$, the equalities 
$(f\circ g,h)= (f,g\circ h)=(g,h\circ f)$ hold.

\begin{Definition}
For any distribution $F \in \S'(\Xi)$ and any function $f \in \S(\Xi)$ we
define
\begin{equation*}
(F\circ f,h):=(F,f\circ h), \quad (f\circ F,h):=(F,h \circ f) \quad 
\hbox{ for all } h \in \S(\Xi).
\end{equation*}
\end{Definition}

The expressions $F \circ f$ and $f \circ F$ are {\it a priori} tempered 
distributions. The Moyal algebra is precisely the set of elements of
$\S'(\Xi)$ that preserves regularity by composition.

\begin{Definition}
The Moyal algebra $\M(\Xi)$ is defined by
\begin{equation*}
\M(\Xi):=\big\{F \in \S'(\Xi)\;|\;F\circ f \in \S(\Xi) \hbox{ and }
f \circ F \in \S(\Xi) \hbox{ for all } f \in \S(\Xi)\big\}.
\end{equation*}
\end{Definition} 
For two distributions $F$ and $G$ in $\M(\Xi)$, the Moyal product 
can be extended by
\begin{equation*}
(F\circ G,h):=(F,G\circ h) \quad \hbox{ for all } h \in \S(\Xi).
\end{equation*}

\begin{Remark}\label{Rem0}
{\rm The set $\M(\Xi)$ with this composition law and the complex 
conjugation $F \mapsto F^\circ$ is a unital ${}^*$-algebra. 
Actually, this extension by duality also gives compositions 
$\M(\Xi)\circ\S'(\Xi)\subset\S'(\Xi)$ and $\S'(\Xi)\circ\M(\Xi)
\subset\S'(\Xi)$. One checks plainly that associativity holds for any 
three factors product with two factors belonging to $\M(\Xi)$ and one 
in $\S'(\Xi)$.}
\end{Remark}

An important result \cite[Prop.~23]{MP2} concerning the Moyal algebra is 
that it contains $\CC^\infty_{\hbox{\tiny \rm pol,u}}(\Xi)$, the space of 
infinitely derivable complex functions on $\Xi$ having uniform 
polynomial growth at infinity. Finally let us quote a result linking 
$\M(\Xi)$ with the functional calculus $\Op^A$ \cite[Prop.~21]{MP2}\;: 
For any vector potential $A$ belonging to $C^\infty_{\hbox{\tiny \rm pol}}(X)$, 
$\Op^A$ is an isomorphism of ${^*}$-algebras between
$\M(\Xi)$ and $\L[\S(X)]\cap \L[\S'(X)]$, where $\L[\S(X)]$ and 
$\L[\S'(X)]$ are, respectively, the spaces of linear continuous 
operators on $\S(X)$ and $\S'(X)$. 

\begin{Remark}\label{Rem}
{\rm We note for further use that very often it is easier to work with 
regularized expressions. For instance, if $f$ and $g$ belong to 
$C^\infty_{\hbox{\tiny \rm pol,u}}(\Xi)$, we can interpret $f\circ g$ as 
the limit $\lim_{m,n\rightarrow\infty}(\chi_n f)\circ(\chi_m g)$, where 
$\chi\in C^\infty_c(\Xi)$ with $\chi(0)=1$ and $\chi_n(\xi):=\chi(\xi/n)$. 
Then $\chi_n f$ is a sequence approximating $f$ in $\S'(\Xi)$ (for example) 
and $(\chi_n f)\circ(\chi_m g)$ is given by the explicit formula 
\eqref{produitMoyal} of the composition law.}
\end{Remark}

\subsubsection*{Twisted crossed product algebras}

Now we recall the definitions of {\it magnetic} twisted 
$C^*$-dynamical systems, of the cor\-responding twisted $C^*$-algebras, and the 
construction of some of their representations in the Hilbert space $\H$. 
These algebras are particular instances of 
more general twisted $C^*$-algebras extensively studied in \cite{BS},
\cite{PR1} and \cite{PR2} (see also references therein). 

For this purpose, let $\A$ be a unital $C^*$-subalgebra of $\BC_u(X)$. 
We shall always assume that $\A$ contains the ideal $\CC_0(X)$ 
and is stable by translations, {\it i.e.}~$\theta_x(a):=a(\cdot +x) \in \A$ 
for all $a \in \A$ and $x \in X$. In the references cited above and in \cite{MPR} 
$\A$ was also assumed to be separable but this is not needed for our developments. 
This algebra can be thought of as a way to encode the anisotropic behaviour of 
the magnetic fields and of the scalar potentials. Thus we consider a magnetic 
field $B$ on $X$ whose components $B_{jk}$ belong to $\A$. 
The expression $\oB$ defined in \eqref{defdeomega} has then some special
properties: For fixed $x$ and $y$, the function $\oB(\cdot;x,y)\equiv \oB(x,y)$ 
belongs to the unitary group $\U(\A)$ of $\A$.
Moreover, the mapping  $X \times X \ni (x,y) \mapsto \oB(x,y) \in \U(\A)$ is a 
2-cocycle on $X$ with values in $\U(\A)$.

The quadruplet $(\A,\theta,\oB,X)$ is a {\it magnetic} example of {\it an abelian 
twisted $C^*$-dynamical system} $(\A,\theta,\omega,X)$. In the general case  
$X$ is an abelian second countable locally compact group, $\A$ is an abelian 
$C^*$-algebra, $\theta$ is a continuous morphism from $X$  to the group of
automorphisms of $\A$,
and $\omega$ is a strictly continuous 2-cocycle with values in the unitary
group of the multiplier algebra of $\A$.
We refer to \cite[Def.~2.1]{MPR} for more explanations.

Given any abelian twisted $C^*$-dynamical system, a natural $C^*$-algebra can 
be defined. We recall its construction.
Let $L^1(X; \A)$ be the set of Bochner integrable functions on $X$ with values
in $\A$, with the $L^1$-norm $\|\phi \|_1 :=\int_X \de x\;\!\|\phi(x)\|_\A$. 
For any $\phi, \psi \in L^1(X;\A)$ and $x \in X$, we define the product
\begin{equation*}
(\phi \diamond \psi)(x):=\int_X\de y\;\theta_{\frac{y-x}{2}}\!\left[\phi(y)
\right]\;\!\theta_{\frac{y}{2}}\!\left[\psi(x-y)\right]\;\!\theta_{-\frac{x}{2}}
\!\left[\omega(y,x-y)\right]
\end{equation*}
and the involution
\begin{equation*}
\phi^{\diamond}(x):=\theta_{-\frac{x}{2}}[\omega(x,-x)^{-1}]\phi(-x)^*.
\end{equation*}
Note that in the magnetic case $\oB(x,-x)=1$. 

\begin{Definition}\label{primel}
The enveloping $C^*$-algebra of $L^1(X,\A)$
is called {\rm the twisted crossed product} and is denoted by
$\A\!\rtimes^\omega_\theta\!X$.
\end{Definition}
 
Let us now consider a continuous vector potential $A$ that generates the 
magnetic field, {\it i.e.}~$A$ is a continuous 1-form on $X$ that satisfies 
$\partial_j A_k - \partial_k A_j = B_{jk}$. 
The relation between $\lambda^A$ defined in equation \eqref{defdelambda} 
and $\oB$ reads (by Stokes Theorem)
\begin{equation}\label{psu}
\lambda^A(q;x)\;\lambda^A(q+x;y)\;\big[\lambda^A(q;x+y)\big]^{-1}=\oB(q;x,y).
\end{equation}  
If $\lambda^A$ were a map $X\ni x \mapsto \lambda^A(\cdot;x) \in \U(\A)$, 
this relation would have said that $\oB$ is a 2-coboundary, or equivalently that
$\oB$ is a trivial 2-cocycle. 
But most the the time this map has only image in $\CC(X;\T)$, the set of 
continuous functions on $X$ with values in the complex numbers of modulus $1$.
For that reason, one says that $\lambda^A$ is a {\it pseudo-trivialization} 
of $\oB$.

Based on relation \eqref{psu}, one can construct a faithful and irreductible 
representation of the algebra $\A\!\rtimes^{\oB}_\theta\!\!X$ in $\B(\H)$, 
that we denote by $\Rep^A$. Equivalently, this corresponds to a covariant
representation of the associated abelian twisted $C^*$-dynamical system.
For each $\phi \in L^1(X;\A)$ and $u \in \H$, the representation is given by
\begin{equation}\label{representation}
\big[\Rep^A(\phi)u\big](x) = \int_X \de y\;\!\lambda^A(x;y-x)\;\!
\phi\big(\hbox{$\frac{1}{2}$}(x+y);y-x\big)\;\!u(y).
\end{equation}
Let us mention that the choice of another vector potential generating
the same magnetic field would lead to a unitarily equivalent 
representation of $\A\!\rtimes^{\oB}_\theta\!\!X$ in $\B(\H)$ (gauge covariance).

By formally comparing \eqref{Op} and \eqref{representation}, one sees
that $\Op^A$ and $\Rep^A$ are connected by a partial
Fourier transformation: $\Op^A(f) = \Rep^A[\FF^{-1}(f)]$, with 
\begin{equation*}
[\FF^{-1}(f)](x,y):=\int_{X^\star}\de p\,e^{-ip\cdot y}f(x,p),
\end{equation*}
for all $x,y \in X$ and suitable $f$.
Then obviously the composition laws $\circ$ and $\diamond$ has to be 
intertwined by $\FF$, {\it i.e.}~$f\circ g = \FF [\FF^{-1}(f)\diamond
\FF^{-1}(g)]$, as it can be checked by a direct computation.
The enveloping $C^*$-algebra $\BBA$ of $\FF\big(L^1(X;\A)\big)$, endowed with 
the multiplication $\circ$ and the complex conjugation, is thus isomorphic to 
$\A\!\rtimes^{\oB}_\theta\!\!X$ via the canonical extension of $\FF$. 
Moreover, one has $\Op^A(\BBA)=\Rep^A(\CBA)$, where $\CBA$ denotes 
for shortness the $C^*$-algebra $\A\!\rtimes^{\oB}_\theta\!\!X$.

It might be here the right place to mention that untwisted crossed products 
are particular cases of groupoids. We suspect that by using {\it twisted} 
groupoids one could get more general results, unifying the present framework 
with the approach of \cite{LN} and \cite{N}.

\subsection{Affiliation}\label{subsecaffil}

In this section we start by recalling the meaning of affiliation, borrowed
from \cite{ABG}. This key concept will then be applied to generalised 
Schr\"odinger operators with magnetic fields. 

\begin{Definition}\label{secundel}
{\rm An observable affiliated to a $C^*$-algebra} $\CCC$ 
is a morphism $\Phi: \CC_0(\R) \to \CCC$.
\end{Definition}

If $\H$ is a Hilbert space and $\CCC$ is a $C^*$-subalgebra of $\B(\H)$,
then a self-adjoint operator $H$ in $\H$ defines an observable 
$\Phi_{\!\hbox{\it \tiny H}}$ affiliated to $\CCC$ if and only if 
$\Phi_{\!\hbox{\it \tiny H}}(\eta) := \eta(H)$ belongs to $\CCC$ for all 
$\eta \in \CC_0(\R)$. 
A sufficient condition is that $(H-z)^{-1} \in \CCC$ for some $z \in \C$ 
with $\im z \neq 0$. Thus an observable affiliated to a $C^*$-algebra is the
abstract version of the functional calculus of a self-adjoint operator.

Given a magnetic field $B$ whose components belong to $\A$, a continuous 
vector potential $A$ that generates $B$ and a suitable symbol 
$h: X^\star \to \R$, our aim is to show that the 
$C_0$-functional calculus of the magnetic Schr\"odinger operator 
$h(\Pi^A)$ (which needs to be carefully defined) belongs to the $C^*$-algebra 
$\Op^A(\BBA)
\subset \B(\H)$. 
The proof of such a statement is rather difficult and we shall do it under 
some smoothness conditions on the magnetic field $B$ and on the symbol $h$.
We point out that we prove in fact a stronger result, 
Theorem \ref{thmaff},  that does not depend on the choice of any 
particular vector potential.

\begin{Definition}\label{tertiel}
\begin{itemize}
\item[{\rm (a)}] For $s \in \R$, a function $h \in \CC^\infty(X^\star)$ is 
{\rm a symbol of type $s$} if the following condition is satisfied:
\begin{equation*}
\forall \alpha \in \N^N, \ \exists\;\! c_\alpha>0 \hbox{ \ such that \ } 
|(\partial^\alpha h)(p)| \leq c_\alpha \langle p\rangle^{s-|\alpha|} \hbox{\ 
for all \ } p \in X^\star,
\end{equation*}
where $\langle p \rangle := \sqrt{1 + p^2}$.  
\item[{\rm (b)}] The symbol $h$ is called {\rm elliptic} 
if there exist $R>0$ and $c>0$ such that 
\begin{equation*}
c\;\!\langle p \rangle^s \leq h(p) \hbox{ \ for all\  } p \in X^\star \hbox{ \
  and\ \  }
|p|\geq R.
\end{equation*}
\end{itemize}
\end{Definition}

We denote by $\,S_{\hbox{\rm \tiny el}}^s(X^\star)\,$ the
family of elliptic symbols of type $s$, and set
$S^\infty_{\hbox{\rm \tiny el}}(X^\star):=\cup_{s}\;\!
S_{\hbox{\rm \tiny el}}^s(X^\star)$. Note that all the classes 
$S^s(X^\star)$ are naturally contained in 
$C^\infty_{\hbox{\tiny \rm pol,u}}(\Xi)$, thus in $\M(\Xi)$. 
For any $z \not \in \R$, we also set $r_z:\R \to \C$ by $r_z(\cdot)
:=(\cdot -z)^{-1}$.

We are in a position to state the results about affiliation.

\begin{Theorem}\label{thmaff}
Assume that $B$ is a magnetic field whose components belong to
$\A\cap BC^\infty(X)$.
Then each real $h \in S^\infty_{\hbox{\rm \tiny el}}(X^\star)$ defines an observable
$\Phi^B_{h}$ affiliated to $\BBA$, such that for any $z \not \in \R$ 
one has 
\begin{equation}\label{fili1}
(h-z)\circ\Phi^B_{h}(r_z)=1=\Phi^B_{h}(r_z)\circ (h-z).
\end{equation} 
In fact one even has $\Phi^B_{h}(r_z)\in\FF \big(L^1(X;\A)\big)\subset \S'(\Xi)$, 
so the compositions can be interpreted as $\M(\Xi)\times\S'(\Xi)\rightarrow\S'(\Xi)$ 
and $\S'(\Xi)\times\M(\Xi)\rightarrow\S'(\Xi)$.
\end{Theorem}

We shall now consider a scalar potential $V\in\A$. It is a standard fact that 
$\A$ consists of multipliers of the algebra $\FF \big(L^1(X;\A)\big)$. 
A straightforward reformulation of the arguments in \cite[p.~365--366]{ABG}
allows then to define the observable $\Phi^B_{h,V}:=\Phi^B_{h}+V$. Considering now $h+V\in\mathcal{S}^\prime(\Xi)$ we remark that we can compute the Moyal product $(h+V-z)\circ\Phi^B_{h,V}(r_z)=(h-z)\circ\Phi^B_{h,V}(r_z)+V\circ\Phi^B_{h,V}(r_z)=1$ 
(by the explicit formula of $\Phi^B_{h,V}$ given in \cite{ABG}). 
This leads to the following statement:

\begin{Corollary}\label{coroll}
We are in the framework of Theorem \ref{thmaff}. Let also $V$ be a real 
function in $\A$. Then $\Phi^B_{h,V}$ is an observable affiliated to $\BBA$, 
such that for any $z \not \in \R$ one has 
\begin{equation*}
(h+V-z)\circ\Phi^B_{h,V}(r_z)=1=\Phi^B_{h,V}(r_z)\circ (h+V-z).
\end{equation*} 
\end{Corollary}

These statements are elegant, being abstract, but in applications one also
needs the represented version:

\begin{Corollary}\label{corol}
We are in the framework of Corollary \ref{coroll}. Let $A$ be a continuous 
vector potential that generates $B$. Then $\ \Op^A(h)+V(Q)\ $ defines a
self-adjoint operator $H_{\!h\!}(A,V)$ in $\H$ with domain given by the image 
of the operator $\Op^A\left[(h-z)^{-1}\right]$ (which do not depend on $z\notin\R$). 
This operator is affiliated to $\Op^A(\BBA) = \Rep^A(\CBA)$. 
\end{Corollary}

In \cite{MP1} we have given an affiliation result for $h(p)=|p|^2$ and 
$\A=BC_{\text{u}}(X)$. In this case we only need that the derivatives 
$\partial ^\alpha B_{jk}$ are bounded for $|\alpha|\le 2$.

\subsection{The essential spectrum}\label{subsecess}

We shall give now a description of the essential spectrum of any observable
affiliated to the $C^*$-algebra $\CBA$. For the
generalised magnetic Schr\"odinger operators of Theorem \ref{thmaff}, 
this is expressed in terms of the spectra of so-called
{\it asymptotic operators}. The affiliation criterion and the algebraic 
formalism introduced above play an essential role in the proof of this result; 
see Section \ref{secess}. We start by recalling some definitions in relation
with topological dynamical systems. 

By Gelfand theory, the abelian $C^*$-algebra $\A$ is isomorphic to the 
$C^*$-algebra $\CC_0(\SA)$, where $\SA$ is the spectrum of $\A$. 
Since $\A$ was assumed unital and contains $\CC_0(X)$, $\SA$ is a
compactification of $X$. We shall therefore identify $X$ with a dense open
subset of $\SA$. 
The group law $\theta: X \times X \to X$ extends then to a continuous map
$\tilde{\theta}: X \times \SA \to \SA$, because $\A$ was also assumed to 
be stable under translations.
Thus the complement $\FA$ of $X$ in $\SA$ is closed and invariant; 
it is the space of a compact topological dynamical system.
For any $\var \in \FA$, let us call the set 
$\{\tilde{\theta}(x,\var) \; | \; x \in X\}$ {\it the orbit generated by
$\var$}, and its closure a {\it quasi-orbit}.
Usually there exist many elements of $\FA$ that generate the same quasi-orbit. 
In the sequel, we shall often encounter the 
restriction $a_\iF$ of an element $a \in \A \equiv \CC(\SA)$ to a quasi-orbit 
$F$. Naturally $a_\iF$ is an element of $\CC(F)$, but we shall
show in Section \ref{secess} that this algebra can  be realized as a 
subalgebra of $\BC_u(X)$. By a slight abuse of notation, we shall
identify $a_\iF$ with a function defined on $X$, thus inducing a 
multiplication operator in $\H$.

The calculation of the essential spectrum may be performed at an abstract
level, {\it i.e.}~without using any representation, as shown in Subsection 
\ref{secab}. In the next statement we present for convenience a represented version.

\begin{Theorem}\label{thmess}
Let $B$ be a magnetic field whose components belong to
$\A\cap BC^\infty(X)$ and let $V\in\A$ be a real function.
Assume that $\{F_\nu\}_\nu$ is a covering of $\FA$ by quasi-orbits.
Then for each real $h \in S^\infty_{\hbox{\rm \tiny el}}(X^\star)$ one has
\begin{equation}\label{aista}
\sigma_{\hbox{\rm \tiny ess}}\big[H_{\!h\!}(A,V)\big] =
\overline{\bigcup_\nu \sigma[H_{\!h\!}(A_\nu,V_\nu)]},
\end{equation}
where $A$, $A_\nu$ are continuous vector potentials for $B$, 
$B_\nu\equiv B_{F_\nu}$, and $V_\nu \equiv V_{F_\nu}$.
\end{Theorem}

The operators $H_{\!h\!}(A_\nu,V_\nu)\equiv h\big(\Pi^{A_{\nu}}\big)+V_\nu$ 
are the asymptotic operators mentioned earlier. We shall show in Section
\ref{secess} that these operators are affiliated to faithful representations 
in $\B(\H)$ of quotients of $\CBA$ by corresponding natural ideals. 
All the spectra appearing in \eqref{aista} are only depending on the respective 
magnetic fields, by gauge covariance. This will be strengthened in Subsection 
\ref{secab} in which a manifestly invariant result will be given in an abstract 
framework.

\subsection{A non-propagation result}\label{subsecprop}

We finally describe how the localization results proved in
\cite{AMP} in the case of Schr\"odinger operators without magnetic field 
can be extended to the situation where a magnetic field is present. 
Once again, the algebraic formalism and the affiliation criterion 
introduced above play an essential role in the proofs: 
see Section \ref{secprop}. 
We first introduce the trace on $X$ of a base of neighbourhoods
of an arbitrary quasi-orbit in $\SA$. 

For any quasi-orbit $F$, let $\NN_{\!F}$ be the family of sets of the
form $W = \W\cap X$, where $\W$ is any element of a base
of neighbourhoods of $F$ in $\SA$. We write $\chi_W$ for the characteristic 
function of $W$. 

\begin{Theorem}\label{thmprop}
Let $B$ be a magnetic field whose components belong to
$\A\cap BC^\infty(X)$, let $V$ be a real scalar potential that 
belongs to $\A$ and let $h$ be a real element of 
$S^\infty_{\hbox{\rm \tiny el}}(X^\star)$. 
Assume that $F\subset F_\A$ is a quasi-orbit. 
Let $A$, $A_F$ be continuous vector potentials for $B$ and $B_F$.
If $\eta \in \CC_0(\R)$ with 
$\text{\rm{supp}}\,(\eta) \cap \sigma[H_{\!h\!}(A_F,V_F)] = \emptyset$,
then for any $\e >0$ there exists $W \in \NN_{\!F}$ such that
\begin{equation*}
\big\|\chi_W(Q)\;\!\eta[H_{\!h\!}(A,V)]\big\|\leq \e.
\end{equation*}
In particular, the inequality 
\begin{equation*}
\big\|\chi_W(Q)\;\!e^{-itH_{\!h\!}(A,V)}\;\!\eta[H_{\!h\!}(A,V)]\;\!u\big\|\leq 
\e \|u\|
\end{equation*} 
holds, uniformly in $t \in \R$ and $u \in \H$.
\end{Theorem}

The last statement of this theorem gives a precise meaning to the notion of 
non-propagation. Heuristically, if the spectral support of $u\in \H$ with 
respect to the operator $H_{\!h\!}(A,V)$ does not meet the 
spectrum of the asymptotic operator corresponding to a quasi-orbit, then the 
state $u$ cannot propagate under the evolution given by $e^{-itH_{\!h\!}(A,V)}$
in the direction of this quasi-orbit. 
We refer to the remark on page 1223 of \cite{AMP} 
for physical explanations and interpretations of this result.

\section{Affiliation}\label{secaffiliation}

In this section we derive our affiliation criterion. In Subsection \ref{subsecaf}
we indicate the main steps of the proof of Theorem \ref{thmaff}. Some technical 
details are included in an appendix. Corollary \ref{corol} is obtained in 
Subsection \ref{subseccor}, as a direct consequence of the theorem. 
We assume tacitly all the hypotheses of Theorem \ref{thmaff}.

\subsection{The proof of the affiliation criterion}\label{subsecaf}

The proof of Theorem \ref{thmaff} will be based on the following strategy:
Let $\M$ be an associative algebra with a composition law denoted by $\circ$ 
and let $\h$ be an element of $\M$.
Our aim is to find the inverse for $\h$.
Assume that $\h'$ is another element such that $\h \circ \h'$ and
$\h' \circ \h$ are invertible. 
These inverses are written $(\h \circ \h')^{(-1)}$ and $(\h' \circ \h)^{(-1)}$
respectively.
Then, the element $\h' \circ (\h \circ \h')^{(-1)}$ is obviously a right
inverse for $\h$ and the element $(\h' \circ \h)^{(-1)} \circ \h'$  a left 
inverse for $\h$. 
Both expressions are thus equal to $\h^{(-1)}$.

In the sequel, we shall take for $\h$ the strictly positive symbol $h+a$, 
with $a$ large enough, and for $\h'$ 
its pointwise inverse $(h+a)^{-1}$. Finding an inverse $(h+a)^{(-1)}$ for $h+a$ 
with respect to the composition law $\circ$ will lead rather easily to 
an observable. In the calculations below we shall use tacitly the approximation 
procedure described in Remark \ref{Rem}. For several arguments we will be forced 
to get out of the algebra $\M=\M(\Xi)$. This will be easily dealt with, 
by a suitable use of elements of $\S'(\Xi)$.

\begin{proof}[Proof of Theorem \ref{thmaff}]

(i) Let us consider an elliptic symbol $h$ of order $s$ and fix some real 
number $a\geq -\inf h +1$. 
We set $h_a:=h+a$, and denote by $h_a^{-1}$ its inverse with respect to
pointwise multiplication, 
{\it i.e.}~$h_a^{-1}(p):=(h(p)+a)^{-1}$ for all $p \in X^\star$.
It is clear that $h_a^{-1}$ is a symbol of type~$-s$.
Since both functions $h_a$ and $h_a^{-1}$ belong to 
$\CC^\infty_{\hbox{\tiny \rm pol,u}}(\Xi)$, and thus to the Moyal algebra 
$\M(\Xi)$, one can calculate their product. 
By using \eqref{produitMoyal} we obtain
\begin{equation}\label{for1}
\left(h_a \circ h_a^{-1}\right)(q,p) = 4^N \int_{X} \de x \int_{X^\star} 
\de k \int_{X} \de y \int_{X^\star} \de l\ \! e^{-2i(k\cdot y-l\cdot x)}
\gB(q;2x,2y)\frac{h_a(p-k)}{h_a(p-l)},
\end{equation}
with $\gB(q;2x,2y):=\oB\big(q-x-y;2x,2(y-x)\big)$. 
The last factor in the integral does not depend on $x$ and $y$; 
it can be developed:
\begin{equation}\label{for2}
\frac{h_a(p-k)}{h_a(p-l)} = 1 + 
\sum_{j=1}^N(l_j-k_j)
\frac{\int_0^1 \de t\;\!(\partial_j h)\big(p-l+t(l-k)\big)}{h(p-l)+a}
=: 1 + \sum_{j=1}^N F_{a,j}(p;k,l)\;.
\end{equation} 
Moreover, let $\gtB(q;k,l) \equiv (\FFF \;\!\gB) (q;k,l)
:=\int_{X} \de x \int_{X} \de y\;\! e^{-ik\cdot y}
\;\!e^{il\cdot x} \;\!\gB(q;x,y)$.
Then the following equality holds (in the sense of distributions, by using 
Remark \ref{Rem}):
\begin{equation}\label{for3}
\int_{X^\star} \de k \int_{X^\star} \de l\;\! \gtB(q;k,l)
=\gB(q;0,0)=1.
\end{equation} 
Thus, by inserting \eqref{for2} and \eqref{for3} into \eqref{for1}, we 
obtain 
\begin{equation*}
h_a \circ h_a^{-1} = 1 + \sum_{j=1}^N f_{a,j},
\end{equation*}
with 
\begin{equation}\label{defdef}
f_{a,j}(q;p):= \int_{X^\star} \de k \int_{X^\star} \de l\;
\gtB(q;k,l)\;\!F_{a,j}(p;k,l) =
\big\langle (\FFF\;\! \gB)(q;\cdot,\cdot),F_{a,j}(p;\cdot,\cdot)\big
\rangle.
\end{equation}
The last notation is used in order to emphasize the duality between 
$\CC^\infty_{\hbox{\rm \tiny pol,u}}(X^\star\times X^\star)$ and its
dual. Indeed, for $q,p$ fixed, Lemma \ref{lemsurF} proves that
$F_{a,j}(p;\cdot,\cdot) \in \CC^\infty_{\hbox{\rm \tiny pol,u}}
(X^\star\times X^\star)$, and Lemma \ref{lemsurgamma} proves that 
$\gB(q,\cdot, \cdot) \in \CC^\infty_{\hbox{\rm \tiny pol}}(X\times X)$, 
so that $(\FFF\;\! \gB)(q;\cdot,\cdot) \in 
\big[\CC^\infty_{\hbox{\rm \tiny pol,u}}(X^\star\times X^\star)\big]'$ 
\cite[Chap.~VII, Thm.~XV]{S}. 

(ii) We are now going to deduce some useful estimates on $f_{a,j}$.
We set $\langle P_x \rangle\equiv\langle -i\partial_x \rangle$.
For $\alpha, j$ fixed and $m,n$ integers that we shall choose below, one has
\begin{equation*}
|(\partial^\alpha_p f_{a,j})(q;p)| \leq
\sup_{x,y \in X}|\langle x \rangle^{-n} \langle y \rangle^{-n}
\langle P_x \rangle^m \langle P_y \rangle^m \;\!\gB(q;x,y)|\  \cdot
\end{equation*}
\begin{equation}\label{jok}
\big\|\langle x \rangle^{-N}\;\!\langle y \rangle^{-N}
\big\|_{L^2(X \times X)} \
\big\| \langle P_k \rangle^{n+N} \langle P_l \rangle^{n+N}
\langle k \rangle^{-m} \langle l \rangle^{-m}\;\!
\left(\partial^\alpha_p F_{a,j}\right)(p;\cdot,\cdot) \big\|_{L^2(X^\star \times
  X^\star)}\;.
\end{equation}
By taking into account \eqref{dom}, subject of Lemma \ref{lemsurF}, and by some 
simple calculations, one can fix $m$ such that the last 
factor of \eqref{jok} is dominated by $c_n\;\! a^{-1/\mu}\;\!\langle p
\rangle^{s/\mu-1-|\alpha|}$, with $\mu >\max\{1,s\}$.
Then, by using Lem\-ma~\ref{lemsurgamma}, one can choose $n$ (depending on $m$) 
such that the first factor on  the r.h.s.~term of \eqref{jok} is bounded. 
Altogether, one obtains
\begin{equation}\label{majorf}
|(\partial^\alpha_p f_{a,j})(q;p)| \leq c\;\! a^{-1/\mu}\;\!\langle p
\rangle^{s/\mu-1-|\alpha|},
\end{equation}
where $c$ depends on $\alpha$ and $j$ but not on $p,q$ or $a$.

(iii) Let us now show that for each $j$, $\FF^{-1}(f_{a,j})$ is an element of 
$L^1(X; \A)$, and thus belongs to the $C^*$-algebra $\CBA$. 
The partial Fourier transform $\FF$ was defined at the 
end of Subsection \ref{fram}.

By taking into account Lemma \ref{lemsurgamma}, the r.h.s.~of the equation
\eqref{defdef} can be rewritten as
$\big \langle \gB(q; \cdot, \cdot), (\FFF^* F_{a,j})(p,\cdot,\cdot)
\big\rangle$, the duality between 
$\CC^\infty_{\hbox{\rm \tiny pol}} (X \times X)$ and 
$\big(\CC^\infty_{\hbox{\rm \tiny pol}}(X\times X)\big)'=\FFF^*C^\infty_{\hbox{
\tiny \rm pol,u}}(X^\star\times X^\star)$. 
As $\gB$ defines a function from $X\times X$ to $\A$ (see Lemma \ref{lemsurgamma}) 
that is of class $C^\infty_{\hbox{\rm \tiny pol}}(X\times X)$, we can easily prove 
that $f_{a,j}(\cdot;p)$ belongs to $\A$, for all $p \in X^\star$ (by using partitions 
of unity on $X\times X$ and by approximating the duality pairing with finite linear 
combinations of elements in $\A$). 

This observation together with \eqref{majorf} imply that
the hypotheses of Lemma \ref{lemdeRadu} are fulfilled for each $f_{a,j}$, 
with $t=-\left(1-s/\mu\right)<0$.
It follows that $\FF^{-1}(f_{a,j})$ belongs to $L^1(X;\A)$ and that there
exists $C>0$ such that
\begin{equation*}
\|\FF^{-1}(f_{a,j})\|_{1}\leq C\;\!a^{-1/\mu}.
\end{equation*}
Thus, for $a$ large enough, the strict inequality $\big\|\sum_{j=1}^N \FF^{-1}(f_{a,j})
\big\|_{1}<1$ holds. It follows that $\FF^{-1}(1+\sum_{j=1}^N f_{a,j})$ is
invertible in $\widetilde{L^1}$, the minimal unitization of $L^1(X;\A)$.
Equivalently, $h_a \circ h_a^{-1} \equiv 1+\sum_{j=1}^N f_{a,j}$ is invertible 
in $\widetilde{\FF (L^1)}$, the minimal unitization of $\FF\big(L^1(X;\A)\big)$. 
Its inverse will be denoted by $\left(h_a \circ h_a^{-1}\right)^{(-1)}$.

(iv) We recall that $h_a^{-1}\in S^{-s}(X^\star)$. Then, by Lemma \ref{lemdeRadu} we get 
that $h_a^{-1}\in \FF\big(L^1(X)\big)\subset\FF\big(L^1(X;\A)\big)$. 
Thus $h_a^{-1}\circ(h_a\circ h_a^{-1})^{(-1)}$ is a well defined element 
of $\FF\big(L^1(X;\A)\big)$.
Moreover, one readily gets $h_a\circ [h_a^{-1}\circ(h_a\circ h_a^{-1})^{(-1)}]=1$. 
For this, just think of $h_a$ and $h_a^{-1}$ as elements 
of the Moyal algebra $\M(\Xi)$ and interpret $(h_a\circ h_a^{-1})^{(-1)}
\in\widetilde{\FF (L^1)}$ as an element of $\S'(\Xi)$. The needed associativity 
follows easily from the definition by duality of the composition law as stated in 
Remark \ref{Rem0}. 
In the same way one obtains $[(h_a^{-1}\circ h_a)^{(-1)}\circ h_a^{-1}]\circ h_a=1$ 
in $\M(\Xi)$. 
In conclusion, there exists $a_0\geq -\inf h+1$ such that for any $a> a_0$ the symbol 
$h_a$ possess an inverse with respect to the Moyal product
\begin{equation*}
h_a^{(-1)}:=h_a^{-1}\circ(h_a\circ h_a^{-1})^{(-1)}=(h_a^{-1}\circ h_a)^{(-1)}\circ 
h_a^{-1}\in\S'(\Xi)
\end{equation*}
that also belongs to $\FF\big(L^1(X;\A)\big)\subset\BBA$. The second equality follows 
from Remark \ref{Rem0} or Remark \ref{Rem} by straightforward arguments.

(v) We define $\Phi^B_h(r_x):=h^{(-1)}_{-x}$ for $x<-a_0$. Then $\Phi^B_h(r_x)\in 
\FF\big(L^1(X;\A)\big) \subset \BBA \cap\mathcal{S}^\prime(\Xi)$, its norm is 
uniformly bounded for $x$ in the given domain and $(h-x)\circ\Phi^B_h(r_x)=
\Phi^B_h(r_x)\circ(h-x)=1$, as shown above. This allows us to obtain an extension 
to the half-strip $\{z=x+iy\mid x<-a_0, |y|<\delta\}$ for some $\delta>0$ by setting
\begin{equation}\label{def-rez}
\Phi^B_h(r_z):=\Phi^B_h(r_x)\circ\{1+(x-z)\Phi^B_h(r_x)\}^{(-1)}.
\end{equation}
It follows that
\begin{equation*}
(h-z)\circ\Phi^B_h(r_z)
=\{(h-x)\circ\Phi^B_h(r_x)+(x-z)\Phi^B_h(r_x)\}\circ\{1+(x-z)\Phi^B_h(r_x)\}^{(-1)}=1.
\end{equation*}
We now prove that the map 
\begin{equation*}
\{z=x+iy\mid x<-a_0, |y|<\delta\}\ni z\mapsto \Phi^B_h(r_z)
\in\FF\big(L^1(X;\A)\big)
\end{equation*} 
satisfies the resolvent equation. Let us choose two complex numbers $z$ and $z'$ 
in this domain and subtract the two equations
\begin{equation}\label{inv}
(h-z)\circ\Phi^B_h(r_z)=1,\qquad (h-z^\prime)\circ\Phi^B_h(r_{z^\prime})=1
\end{equation}
in order to get $(h-z)\circ\{\Phi^B_h(r_z)-\Phi^B_h(r_{z^\prime})\}+(z^\prime-
z)\Phi^B_h(r_{z^\prime}) =0$.
By multiplying at the left with $\Phi^B_h(r_z)$ and by using the associativity, 
we obtain the resolvent equation
\begin{equation*}
\Phi^B_h(r_z)-\Phi^B_h(r_{z^\prime})=(z-z^\prime)\Phi^B_h(r_z)\circ
\Phi^B_h(r_{z^\prime}).
\end{equation*}
Now, setting $z^\prime=\overline{z}=x-iy$ with $y>0$ and taking norms we get
\begin{equation*}
\|\Phi^B_h(r_z)\|_{\FF(L^1(X;\A))}\leq y^{-1}. 
\end{equation*}
With this estimate and formula \eqref{def-rez}, the function 
$z\mapsto\Phi^B_h(r_z)$ can be extended to the domain $\C\setminus [-a_0,+\infty)$, 
preserving the relations \eqref{inv}.
The resolvent equation may be proved in a similar way to hold on the entire domain 
$\C\setminus [-a_0,+\infty)$ and analyticity of the defined function follows in an 
evident way.

(vi) Thus we have got an analytic map $\C\setminus [-a_0,+\infty)\ni z\rightarrow\Phi^B_h(r_z)
\in\FF\big(L^1(X;\A)\big)$ satisfying the resolvent equation and the symmetry condition. 
A general argument presented in \cite[p.~364]{ABG} allows now to extend in a unique way 
the map $\Phi^B_h$ to a $C^*$-algebra morphism $C_0(\R)\rightarrow\BBA$.
\end{proof}

\subsection{The represented version}\label{subseccor}

This subsection consists only in the proof of the represented version on Theorem 
\ref{thmaff}.
\begin{proof}[Proof of Corollary \ref{corol}]
We shall first consider the case $V=0$ and then add $V$ as a bounded perturbation.

Let us denote by $\D_z$ the range of the operator $\Op^A[\Phi^B_{h}(r_z)]\in \B(\H)$. 
By the resolvent identity it follows immediately that it is a subspace of $\H$ that 
does not depend on $z\in\C \setminus \R$. Thus we set $\D_z\equiv\D$.
Since $h\in\M(\Xi)$, one has $\Op^A(h)\in\mathcal{L}[\mathcal{S}(X)]\cap
\mathcal{L}[\mathcal{S^\prime}(X)]$. We interpret it as a linear operator in 
$\mathcal{S^\prime}(X)$ and set $H_{\!h\!}(A,0):=\Op^A(h)|_\D$.

Now, by applying $\Op^A$ to \eqref{fili1} we get 
\begin{equation*}
\{H_{\!h\!}(A,0)-z\boldsymbol{1}\}\Op^A[\Phi^B_{h}(r_z)]=\boldsymbol{1}
\end{equation*}
and
\begin{equation*}
\Op^A[\Phi^B_{h}(r_z)]\{\Op^A(h)-z\boldsymbol{1}_{\mathcal{S}(X)}\}=
\boldsymbol{1}_{\mathcal{S}(X)}.
\end{equation*}
The first identity shows that $H_{\!h\!}(A,0)\D\subset\H$. Straightforwardly 
it is hermitian. 
The second equality implies that $\mathcal{S}(X)\subset\D$ and thus $\D$ is 
dense in $\H$. By the first equality above the ranges of $H_{\!h\!}(A,0)\pm i$ 
both coincide with $\H$. Thus, by the fundamental criterion of self-adjointness, 
$H_{\!h\!}(A,0)$ is self-adjoint.

By construction, $\{\Op^A[\Phi^B_{h}(r_z)]\,|\,z\in \C\setminus \R\}$ is 
the resolvent family of $H_{\!h\!}(A,0)$, which is therefore affiliated to 
$\Op^A(\BBA)$.

Then we define the standard operator sum $H_{\!h\!}(A,V):=H_{\!h\!}(A,0)+V:\D\rightarrow\H$. 
Using the second resolvent equation and the Neumann series the conclusion of the Corollary 
follows easily using \cite[Prop.~2.6]{MPR} as in \cite{MP1}. A different proof could start 
from the result of Corollary \ref{coroll}.
\end{proof}

\section{The essential spectrum}\label{secess}

In this section, we shall consider certain abelian twisted $C^*$-dynamical 
system $(\A, \theta, \omega, X)$ and explain how to calculate the essential 
spectrum of any observable affiliated to the twisted crossed product algebra $\A\!\times^\omega_\theta\!X$. This result is contained in
Proposition \ref{Propess}. 
Then, by using the concrete affiliation criterion obtained in Subsection \ref{subsecaffil}, 
we shall particularize the result to the case of magnetic Schr\"odinger operators 
and prove Theorem \ref{thmess}.    

We start by recalling some definitions in relation with spectral analysis in a
$C^*$-algebraic framework. 
Let $\pi : \CCC \to \CCC'$ be a morphism between two $C^*$-algebras and 
$\Phi$  an observable affiliated to $\CCC$. 
Then $\pi[\Phi]: \CC_0(\R) \to \CCC'$ given by 
$\big(\pi[\Phi]\big)(\eta):=\pi[\Phi(\eta)]$ is an observable affiliated to
$\CCC'$, called {\it the image of $\Phi$ through $\pi$}. 
If $\KK$ is an ideal of $\CCC$, {\it the $\KK$-essential spectrum of $\Phi$} is
\begin{equation*}
\sigma_\KK(\Phi):=
\big\{\lambda \in \R\;|\; \hbox{ if }\eta \in \CC_0(\R) \hbox{ and } 
\eta(\lambda) \neq 0, \hbox{ then }\Phi(\eta) \not \in \KK\big\}.
\end{equation*} 
If $\pi$ denotes the canonical morphism $\CCC \to \CCC/\KK$, one has 
$\sigma_\KK(\Phi) = \sigma_{\{0\}}(\pi[\Phi])$.

In the particular situation when $\CCC$ is a $C^*$-subalgebra of $\B(\H)$ 
for some Hilbert space $\H$, any self-adjoint operator $H$ in $\H$ 
defines an observable $\Phi_{\!\hbox{\it \tiny H}}$ affiliated to $\CCC$ 
by its functional calculus $C_0(\R)\ni\eta\mapsto \eta(H)
\equiv\Phi_{\!\hbox{\it \tiny H}}(\eta)$ if and only if $\Phi_{\!\hbox{\it \tiny H}}
(r_z)\in\CCC$ for some $z \notin \R$. Then $\sigma_{\{0\}}(\Phi_{\!\hbox{\it \tiny H}})$ 
is the usual spectrum $\sigma(H)$ of $H$. 
Moreover, if $\CCC$ contains the ideal $\K(\H)$ of compact operators on
$\H$, then $\sigma_{\K(\H)}(\Phi_{\!\hbox{\it \tiny H}})$ is equal to the 
essential spectrum $\sigma_{\hbox{\rm \tiny ess}}(H)$ of $H$. 
Here we shall be mainly interested in the usual spectrum and in the essential 
spectrum. The need for the $\KK$-essential spectrum with $\KK$ different from $\{0\}$ or 
$\K(\H)$ will appear only in Section \ref{secprop}.

\subsection{The abstract construction}\label{secab}

In this subsection $(\A,\theta,\omega,X)$ will be an abelian twisted $C^*$-dynamical 
system. Thus $X$ is an abelian, second countable locally compact group and $\A$ an abelian, 
unital $C^*$-subalgebra of $\BC_u(X)$ stable under translations and containing 
$C_0(X)$. We recall that the spectrum $\SA$ of $\A$ is a
compactification of $X$, endowed with an action $\tilde{\theta}$ of $X$ by 
homeomorphisms. For any quasi-orbit $F$ we define
$$\A^F:=\big\{a \in \CC(\SA)\ |\ a|_{F}=0\big\}.$$
By identifying $\A$ with $\CC(\SA)$, $\A^F$ will be an invariant
ideal of $\A$.
Obviously the unitary group $\U(\A^F)$ of the multiplier algebra
of $\A^F$ contains the unitary group $\U(\A)$ of $\A$.
Consequently, the abelian twisted dynamical system $(\A^F, \theta,\omega, X)$
obtained by replacing $\A$ with $\A^F$ and performing suitable restrictions is well 
defined. Furthermore, the twisted crossed product
$\A^F\!\!\rtimes^{\omega}_\theta\!X$
may be identified with an ideal of $\A\!\rtimes^{\omega}_\theta\!X$
\cite[Prop.~2.2]{PR2}.

In order to have an explicit description of the quotient,
let us first note that $\A / \A^F$ is canonically isomorphic to the unital
$C^*$-algebra $\CC(F)$ of all continuous functions on $F$.
The natural action of $X$ on $b \in \CC(F)$ is given by $(\theta_x b)(\var) =
b[\tilde{\theta}(x,\var)]$ for each $x \in X$ and $\var \in F$.
Now, for each $x,y \in X$, the restriction of $\omega(x,y) \in \U(\A)$ to $F$
gives rise to a 2-cocycle $\omega_{\!\iF}: X \times X \to \U\big(\CC(F)\big)$.
Thus $\big(\CC(F), \theta, \omega_{\!\iF},X\big)$ is a well-defined abelian twisted 
$C^*$-dynamical system. Moreover the quotient $\A\!\rtimes^{\omega}_\theta\!X /
\A^F\!\!\rtimes^{\omega}_\theta\!X$ may be
identified with the corresponding twisted crossed product
$\CC(F)\!\rtimes^{\omega_{\!F}}_\theta\!X$. This follows from \cite[Prop.~2.2]{PR2} 
if $\A$ is separable. For the non-separable case, just perform obvious modifications 
in the proof of \cite[Th.~2.10]{GI3} to accommodate the $2$-cocycle.
Let us recall that $a_{\!\iF}$ denotes the restriction of 
$a \in \A\equiv \CC(\SA)$ to $F$. 
Then the image of $\phi \in L^1(X;\A)$ through the canonical morphism
$\pF:\A\!\rtimes^{\omega}_\theta\!X  \to
\CC(F)\!\rtimes^{\omega_{\!\iF}}_\theta\!X$ is the element of
$L^1\big(X;\CC(F)\big)$ given by $\big(\pF[\phi]\big)(x) = [\phi(x)]_{\!\iF}$
for all $x \in X$.

Let us consider a covering $\{F_{\!\nu}\}_\nu$ of $\FA$ by quasi-orbits.
At the algebraic level, the covering requirement reads
$\cap_\nu \A^{F_{\!\nu}}=\CC_0(X)$.
It implies the equality
\begin{equation*}
\bigcap_\nu \big(\A^{F_{\!\nu}}\!\!\rtimes^{\omega}_\theta
\!X\big) = \CC_0(X)\!\rtimes^{\omega}_\theta\!X\;.
\end{equation*}
By putting all these together one obtains, {\it
cf.}~\cite[Prop.~1.5]{M1}:

\begin{Proposition}\label{Propess}
Let $\{F_{\!\nu}\}_\nu$ be a covering of $\FA$ by quasi-orbits.
\begin{enumerate}
\item[\rm{(i)}] There exists an injective morphism
\begin{equation*}
\A\!\rtimes^{\omega}_\theta\!X \ / \
\CC_0(X)\!\rtimes^{\omega}_\theta\!X \hookrightarrow
\prod_\nu \CC(F_{\!\nu})\!\rtimes^{\omega_{\!\iF_{\!\nu}}}_\theta\!X\;.
\end{equation*}
\item[\rm{(ii)}] If $\Phi$ is an observable affiliated to
$\A\!\rtimes^{\omega}_\theta\!X$ and $\pi_{\!\iF_{\!\nu}}$ denotes the 
canonical surjective morphism $\A\!\rtimes^{\omega}_\theta\!X \to
\CC(F_{\!\nu})\!\rtimes^{\omega_{\!\iF_{\!\nu}}}_\theta\!X$, then, with $\KK:= \CC_0(X)\!\rtimes^{\omega}_\theta\!X$, we have
\begin{equation}\label{spectre}
\sigma_{\KK}(\Phi) =
\overline{\bigcup_\nu \sigma(\pi_{\!\iF_{\!\nu}}[\Phi])}\;.
\end{equation}
\end{enumerate}
\end{Proposition}

We now introduce a represented version of this proposition in the
Hilbert space $\H$.
Let $\lambda \in C\big(X;\CC(X;\T)\big)$ be a 1-cochain satisfying the 
relation 
\begin{equation}\label{1cochain}
\lambda(x)\;\!\theta_x[\lambda(y)]\;\!\lambda(x+y)^{-1}=\omega(x,y) \quad 
\hbox{ for all } x,y \in X.
\end{equation}
It was proved in \cite[Prop~2.14]{MPR} that such a pseudo-trivialization
function $\lambda$ always exists.
The associated representation of $\A\!\rtimes^\omega_\theta\!X$ in $\B(\H)$
defined by \eqref{representation}, but with $\lambda^A$ replaced by $\lambda$, is denoted by $\Rep^\lambda$.
We recall from \cite[Prop.~2.17]{MPR} that $\Rep^\lambda$ is irreducible 
and faithful and that
$\Rep^\lambda\big(\CC_0(X)\!\rtimes^{\omega}_\theta\!X\big)$ 
is equal to $\K(\H)$.
If $\Phi$ is an observable affiliated to $\A\!\rtimes^\omega_\theta\!X$, then 
the l.h.s.~term of \eqref{spectre} is equal to
$\sigma_{\hbox{\rm \tiny ess}}\big(\Rep^\lambda(\Phi)\big)$, and it does not depend on a particular choice of $\lambda$.

In order to construct a faithful representation of
$\CC(F_{\!\nu})\!\rtimes^{\omega_{\!\iF_{\!\nu}}}_\theta\!\!X$ in $\H$,
we rely on the natural realization of the restriction of $\A$ 
to a quasi-orbit mentioned in Subsection \ref{subsecess}.  
Let $F$ be a quasi-orbit and $\var$ an element of $\FA$ that generates it. 
Then, for any $b \in \CC(F)$ and $x \in X$, set $b_{\var}(x):=
b[\tilde{\theta}(x,\var)]$.
By taking into account the surjectivity of the morphism $\A \to \CC(F)$ and
the continuity of translations in $\A \subset \BC_u(X)$, one easily sees that
$b_{\var}: X \to \C$ belongs to $\BC_u(X)$.
Furthermore, the induced action of $X$ on $b_{\var}$ coincides with the natural
action of $X$ on $\BC_u(X)$.
One has thus obtained an embedding of $\CC(F)$ in $\BC_u(X)$.
By an abuse of notation, we shall keep writing $b$ for $b_\var$, and $C(F)$
for the corresponding $C^*$-subalgebra of 
$\BC_u(X)$.

Now, by choosing any 1-cochain $\lambda^F \in C\big(X;\CC(X;\T)\big)$ 
satisfying the pseudo-triviality relation \eqref{1cochain} with 
$\lambda = \lambda^F$ and $\omega = \omega_{\!\iF}$,
one can construct the faithful Schr\"odinger representation $\Rep^{\lambda^F}$ of the
algebra $\CC(F)\!\rtimes^{\omega_{\!F}}_\theta\!X$.
Thus, if $\Phi$ is the observable affiliated to $\A\rtimes^\omega_\theta\!X$
of Proposition \ref{Propess}, then each observable 
$\pi_{\!\iF_{\!\nu}}(\Phi)$ can be represented as a observable affiliated 
to a $C^*$-subalgebra of $\B(\H)$ and having the same spectrum.
This remark makes the calculation of the r.h.s.~terms in
\eqref{spectre} more concrete. The particular case treated in Theorem~\ref{thmess} 
is proved now.

\subsection{Application to magnetic Schr\"odinger operators}

We particularize the above construction to the case of a magnetic 2-cocycle $\oB$.
So, we consider a magnetic field $B$ whose components belong to $\A$.
We shall need the following parametrized formula:
for $q,x,y \in X$ 
\begin{equation}\label{defdeux}
\oB(q;x,y)  =\exp\Big\{-i\sum_{j,k=1}^Nx_j\;\!y_k
\int_0^1 \de s \int_0^1 \de t\;\!s\;\!B_{jk}(q + s\;\!x +s\;\!t\;\!y)\Big\}.
\end{equation}

We are now in a position to prove Theorem \ref{thmess}. It consists essentially
in an application of Proposition \ref{Propess} together with a partial Fourier
transformation. 

\begin{proof}[Proof of Theorem \ref{thmess}]
Let us fix a quasi-orbit $F_{\!\nu}$; obviously $\oB|_{F_{\!\nu}}=
\omega^{\hbox{\it \tiny B}_{\!\iF_{\!\nu}}}$ with natural identifications. 
Then the morphism 
\begin{equation*}
\FF\big(L^1(X;\A)\big) \ni f \mapsto \FF\big(\pi_{\!\iF_{\!\nu}}[\FF^{-1}(f)]\big)
\in \FF\Big(L^1\big(X;C(F)\big)\Big)
\end{equation*}
extends to a surjective morphism 
$\tilde{\pi}_{\!\iF_{\!\nu}}:\BBA \to {\mathfrak B}^{B_{\!F_{\!\nu}}}_{\!C(F_\nu)}$.
The equality \eqref{spectre} can then be rewritten in the framework of 
$\BBA$ and for the observable $\Phi^B_{h,V}$:
\begin{equation*}
\sigma_{\hbox{\tiny \rm ess}}(\Phi^B_{h,V}) =
\overline{\bigcup_\nu \sigma\big(\tilde{\pi}_{\!\iF_{\!\nu}}[\Phi^B_{h,V}]\big)}.
\end{equation*}
The result follows now from the central observation 
that $\tilde{\pi}_{\!\iF_{\!\nu}}[\Phi^B_{h,V}]$ is equal to 
$\Phi^{B_{\!F_{\!\nu}}}_{h,V_{\!F_{\!\nu}}}$, by considering faithful 
representations (i) of $\BBA$ through $\Op^{A}$ and (ii) of 
${\mathfrak B}^{B_{\!F_{\!\nu}}}_{\!C(F_{\!\nu})}$ through $\Op^{A_\nu}$ and 
by applying Corollary \ref{corol}.
\end{proof}

\section{Non-propagation properties}\label{secprop}

As mentioned earlier, the result of non-propagation is mainly an 
adaptation of \cite{AMP} in the presence of a magnetic field together with 
the use of an approximate unit introduced in \cite{M2}. 
Since all notations and concepts have already been introduced, it only remains
to prove Theorem \ref{thmprop}. We start by recalling an easy result 
of \cite[Lem.~1]{AMP}.

\begin{Lemma}\label{lemAMP}
Let $\KK$ be an ideal in a $C^*$-algebra $\CCC$ and $\Phi$ an observable
affiliated to $\CCC$. If $\eta \in \CC_0(\R)$ and $\eta(\lambda)=0$ for all
$\lambda \in \sigma_\KK(\Phi)$, then $\Phi(\eta) \in \KK$.
\end{Lemma}

\begin{proof}[Proof of Theorem \ref{thmprop}]
Let $\KK := {\mathfrak B}^B_{\!\!\A^F} \equiv \FF
\big({\mathfrak C}^B_{\!\!\A^F}\big)$, 
the ideal of $\BBA$ related to the quasi-orbit $F$, and let
$\tilde{\pi}_{\!\iF}: \BBA \to {\mathfrak B}^{B_{\!F}}_{\!C(F)}$ be the corresponding 
morphism of kernel $\KK$.
We consider the observable $\Phi_{h,V}^B$ that is affiliated
to $\BBA$ by Theorem  \ref{thmaff}.
Then, by taking into account the equality
 $\sigma_\KK(\Phi_{h,V}^B) = 
\sigma\big(\tilde{\pi}_{\!\iF}\big[\Phi_{h,V}^B\big]\big)
= \sigma \big(\Phi_{h,V_{\!F}}^{B_{\!F}}\big)$, 
the hypothesis on $\eta$ and Lemma \ref{lemAMP}, we see that 
$\Phi_{h,V}^B(\eta)$ belongs to $\KK$.

By representing faithfully $\BBA$ in $\B(\H)$ through $\Op^A$ one has
that $\Op^A\big(\Phi_{h,V}^B(\eta)\big)$ belongs to the ideal
$\Op^A(\KK)$.
For the final step of the proof, one only has to remark that the family
$\{\mathbf{1}-\chi_W(Q)\}_{W \in \NN_{\!F}}$ is an approximate unit in 
$\B(\H)$ for $\Op^A\big({\mathfrak B}^B_{\!\!\A^F}\big)\equiv 
\Rep^A\big({\mathfrak C}^B_{\!\!\A^F}\big)$, which is 
straightforward by the description of this type of algebras given in 
\cite[Prop.~2.6]{MPR}.  
\end{proof}

\section{Examples}\label{sexemples}

In this last section, we illustrate Theorem \ref{thmess} on the
essential spectrum by choosing concrete examples of algebras $\A$.
A similar transcription of Theorem \ref{thmprop} on propagation for
these concrete situations could also be performed.
Since an adaptation for the magnetic case of the examples given in
\cite{AMP} is rather straightforward, we leave this to the 
reader.

It is always assumed in the sequel that the components of the magnetic
field $B$ belong to $\A\cap BC^\infty(X)$ and that the scalar potential $V$
belongs to $\A$. It is convenient to write $\sigma[H_{\!h\!}(B,V)]$ for $\sigma[H_{\!h\!}(A,V)]$ and $\sigma_{\hbox{\rm \tiny ess}}
[H_{\!h\!}(B,V)]$ for $\sigma_{\hbox{\rm \tiny ess}}[H_{\!h\!}(A,V)]$ if 
$B=dA$. This is justified by the independence of these sets on a choice of 
a vector potential and, especially, by the abstract approach of Subsection 
\ref{secab}.

The easiest and best known situation is certainly when the algebra $\A$ 
is equal to $\C+C_0(X)$.
In this situation $\A/C_0(X)\cong \C$ and one has $\sigma_{\hbox{\rm \tiny ess}}[H_{\!h\!}(B,V)]=\sigma[H_{\!h\!}(B_\infty, V_\infty)]=
\sigma[H_{\!h\!}(B_\infty,0)]+V_\infty$,
where $B_\infty$, $V_\infty$ are respectively the limits of $B$ and $V$ at infinity. 
For instance, if $h(p)=|p|^2$ (giving the usual magnetic Schr\"odinger operator) 
in $X=\R^2$, we have for $B_\infty\ne0$: $\sigma_{\hbox{\rm \tiny ess}}[H(B,V)]=(2\N+1)B_\infty+V_\infty$, 
a translation by $V_\infty$ of the familiar Landau levels. For  $B_\infty=0$ we clearly obtain 
$\sigma_{\hbox{\rm \tiny ess}}[H(B,V)]=[V_\infty,\infty)$. Some related results may be found in \cite{P}.

We shall now consider more complicated examples.

\subsection{Vanishing oscillation}\label{VOX}

We take $\A$ to be the algebra $VO(X)$ of
{\it vanishing oscillations functions} :

\begin{Definition}
A bounded and uniformly continuous function $a$ belongs to $VO(X)$
if for any $x \in X$, the difference $\theta_x[a]-a$ 
belongs to $C_0(X)$.
\end{Definition}

Obviously, $VO(X)$ is a unital $C^*$-algebra containing $C_0(X)$
and stable by translations.         
It contains also $C^{\hbox{\rm \tiny rad}}(X)$, the algebra of continuous
functions 
that can be extended continuously to the radial compactification of $X$
obtained by adding a sphere at infinity.
But $VO(X)$ is in fact much larger than $C^{\hbox{\rm \tiny rad}}(X)$.
For example, it also contains the set of all bounded $C^1$-functions
with derivatives in $C_0(X)$. 
A simple typical example is $a(x):=f\big((1+|x|)^s\big)$ (suitably regularized 
at the origin), where $f$ is a periodic $C^1$-function of one variable and 
$s$ is a real number strictly smaller than $1$.

To understand what the asymptotic operators should be, let us introduce 
the notion of {\it asymptotic range} of a real, bounded and continuous function 
$\varphi$ defined on $X$. We write $\lambda\in\varphi(X)_{\hbox{\rm \tiny asy}}$
if and only if for any $\e>0$, $\varphi^{-1}[(\lambda-\e,
\lambda+\e)]$ is not relatively compact in $X$. Equivalently, 
$\lambda \in [\lim \inf_{x\rightarrow\infty} \varphi(x),\lim \sup_{x\rightarrow\infty} 
\varphi(x)]$, or it exists a divergent sequence $\mathbf{x}=\{x_n\}_{n \in \N}$ such that $\varphi(x_n)\rightarrow\lambda$ when $n\rightarrow\infty$. 
We recall that a divergent sequence $\{x_n\}_{n \in \N}$ consists in a sequence
of $x_n \in X$ such that $x_n \to \infty$ as $n \to \infty$.
The interest in the set $\varphi(X)_{\hbox{\rm \tiny asy}}$ lies in the fact 
that for any $\A\equiv C(\SA)$ containing $\varphi$, the range of the restriction to 
$\FA$ is exactly $\varphi(X)_{\hbox{\rm \tiny asy}}$, 
{\it i.e.}~$\varphi(\FA)=\varphi(X)_{\hbox{\rm \tiny asy}}$ with a loose notation. 

A nice feature of $VO(X)$ is that it is the largest unital translational
invariant $C^*$-subalgebra of $\BC_u(X)$ such that all quasi-orbits situated at 
infinity are reduced to points. This means that $R\equiv F_{VO(X)}$ admits the 
partition $R=\sqcup_{\var\in R}\{\var\}$ in (quasi-)orbits and $\varphi\mapsto\big(\varphi(\var)\big)_{\var\in R}$ determines the embedding of 
$VO(X)/C_0(X)$ into $\prod_{\var\in R}\C$. 
Using all these in conjunction with Theorem \ref{thmess} leads to 
\begin{equation*}
\sigma_{\hbox{\rm \tiny ess}}[H_{\!h\!}(B,V)]=\overline{\bigcup_{\var\in R}
\sigma[H_{\!h\!}(B(\var), V(\var))]}=\overline{\bigcup_{\mathbf{x}\in \mathcal R}
\sigma[H_{\!h\!}(B_\mathbf{x}, V_\mathbf{x})]},
\end{equation*} 
where the second union is performed over the set $\mathcal R$ of all divergent 
sequences $\mathbf{x} = \{x_n\}_{n\in \N}$ such that there exist a constant magnetic 
field $B_\mathbf{x}$ and a number $V_\mathbf{x}$ satisfying 
$\sup_{j,k}|(B(x_n)-B_\mathbf{x})_{jk}| \to 0$ and $|V(x_n)-V_\mathbf{x}|\to 0$
as $n \to \infty$. We say that $B_\mathbf{x}$ and $V_\mathbf{x}$ are 
{\it asymptotic values} for $B$ and $V$ respectively. 
Various particularizations are available.

\subsection{Comparison with the results of \cite{HM}}

The results of \cite{HM} are very interesting because large classes of unbounded 
potentials and magnetic fields are admitted. In the bounded case, however, they 
are entirely confined to the vanishing oscillation type of anisotropy, as we now 
argue.

For the comparison with the results of \cite{HM}, we need

\begin{Lemma}\label{lemsurVOX}
Let $r \in \N$ and $f \in \BC(X)\cap \CC^r(X)$. 
Assume that $\partial^\alpha f \in \CC_0(X)$ for all $\alpha \in \N^N$
with $|\alpha|=r$. Then $f$ belongs to $VO(X)$ and $\partial^\beta
f \in 
\CC_0(X)$ for all $\beta \in \N^N$  with $1\leq |\beta|\leq r-1$.
\end{Lemma} 

\begin{proof}
Since we were not able to locate this result in the literature, we
sketch its 
proof. Let us first state three remarks which are easily proved.
(i) Under the hypotheses on $f$, one has $\partial^\beta f \in
\BC(X)$ for 
all $\beta \in \N^N$  with $1\leq |\beta|\leq r-1$, {\it cf.}~for
example 
\cite{L}.
(ii) If $g \in \BC^1(X)$ and $\partial_j g \in \CC_0(X)$ for all
$j \in \{1,\dots,N\}$, then $g \in VO(X)$.
(iii) If $h \in \BC^1(X)$ and $\partial_j h \in VO(X)$, then 
$\partial_j h \in \CC_0(X)$.  

Now, if $r=1$, the result is obtained by 
(ii). If $r\geq 2$, let $\beta \in \N^N$ with $|\beta|=r-2$ and set 
$h:=\partial^\beta f \in \BC^2(X)$ by (i). 
For each $k\in \{1,\dots,N\}$, $\partial_k h$ belongs to $VO(X)$ by
(ii), 
and then to $\CC_0(X)$ by (iii). 
By varying $\beta$ and $k$, one obtains $\partial^\gamma f \in \CC_0(X)$
for all $\gamma \in \N^N$ with $|\gamma|=N-1$. 
A bootstrap argument leads to the result.
\end{proof}

We describe now the results of \cite{HM} with slightly modified notations. 
They consider magnetic Schr\"odinger operators $H_{\!h\!}(A,V)$ for the particular 
case $h(p)\equiv h_0(p) :=|p|^2$. The data $V$ and $B=dA$ are subject to the 
following assumptions. There exist $q,r \in \N$ and some smooth function 
$\rho:X\rightarrow [1,\infty)$ with $\rho(x)\rightarrow\infty$ when 
$|x|\rightarrow\infty$, which is also tempered in a sense that is not important here 
\cite[eq.~(1.14)]{HM}, such that
\begin{enumerate}
\item[\rm{(i)}] $V=V_0+\sum_{l=1}^q V^2_l$,
\item[\rm{(ii)}] $V_0\geq -C_1$, $V_0\in C^1(X)$ and $V_l\in C^{r+2}(X)$ for $l=1,\dots,q$,
\item[\rm{(iii)}] $\sum_{|\alpha|=1}|\partial^\alpha V_0|+\sum_{|\alpha|=r+2}\sum_{l=1}^q
|\partial^\alpha V_l|\le C_2\;\varphi^{-1}$,
\item[\rm{(iv)}] for all $j,k$, $B_{jk}\in C^{r+3}(X)$ and $\sum_{|\alpha|=r+1}^{r+3}
\varphi^{|\alpha|-r-1}|\partial^\alpha B_{jk}|\le C_3\;\varphi^{-1}$.
\end{enumerate}
Under these assumptions it is proved that 
\begin{equation}\label{final}
\sigma_{\hbox{\rm \tiny ess}}[H_{\!h_0\!}(B,V)]=\overline{\bigcup_{\mathbf{x}\in\mathcal R}
\sigma[H_{\!h_0\!}(A_\mathbf{x}, V_\mathbf{x})]},
\end{equation}
where $\mathcal R$ is the set of divergent sequences $\mathbf x=\{x_n\}_{n \in \N}$ 
such that the following limits exist:
\begin{enumerate}
\item[\rm{(a)}] $v_0=\underset{n}{\lim}\;V_0(x_n)$ and $v_l^\alpha=\underset{n}{\lim}\;
(\partial^\alpha V_l)(x_n)$ for $l=1,\dots,q$ and $|\alpha|\le r+1$,
\item[\rm{(b)}] $B^\alpha=\underset{n}{\lim}\;(\partial^\alpha B)(x_n)$ for $|\alpha|\le r$.
\end{enumerate}
Then the asymptotic operators $H_{\!h_0\!}(A_\mathbf{x}, V_\mathbf{x})$ are constructed 
with the scalar potential
\begin{equation*}
V_{\mathbf{x}}(x):=v_0+\sum_{l=1}^q\Big(\sum_{|\alpha|\le r+1}\frac{v_l^\alpha}{\alpha!}
x^\alpha\Big)^2
\end{equation*}
and the magnetic potential
\begin{equation*}
A_{\mathbf{x}}(x):=\sum_{|\alpha|\le r}\frac{B^\alpha\cdot x}{\alpha!(2+|\alpha|)}
x^\alpha\;.
\end{equation*}

Let us see how the hypotheses and the conclusion look like when $V$ and $B$ are bounded. 
We ignore the temperedness condition; the fact that $\varphi$ diverges at infinity implies 
that the l.h.s.~of the conditions (iii) and (iv) belong to $C_0(X)$. 
Lemma \ref{lemsurVOX} can be applied and thus $V$ and $B_{jk}$ belong to $VO(X)$. 
In the bounded case the anisotropy covered by \cite{HM} is surely of the vanishing oscillation 
type.

To understand the conclusion under the extra condition that $V$ and $B$ are bounded, note 
that the same Lemma \ref{lemsurVOX} says that all the derivatives of strictly positive order 
are in $C_0(X)$, thus the only non-null constant coefficients in (a) and (b) are those 
corresponding to $\alpha=0$. Then \eqref{final} coincides with our result described above. 

\subsection{Mixed algebras}\label{moreexamples}

In the examples developed above, the quasi-orbits are reduced to singletons. 
We shall introduce some algebras with more complicated  quasi-orbits, leading to 
non-trivial asymptotic operators with variable coefficients. 
There is also a very nice type of anisotropy studied in \cite{GI2} and \cite{GI3} 
under the name {\it potentials belonging to the bumps algebra}. 
It would be interesting to work out the magnetic counterpart. 

We first reconsider an example introduced
in \cite{M1}, to which we refer for details and comments.
Let us introduce the algebra $AP(X)$ of all continuous, 
almost periodic functions on $X$ \cite[16.2.1]{D} :

\begin{Definition}
A bounded and continuous function $a$ on $X$  belongs to $AP(X)$
if and only if it satisfies one of the following equivalent condition :
\begin{itemize}
\item[{\rm (a)}] The set $\{\theta_x[a]\, | \, x \in X\}$ is relatively
compact in $BC(X)$.
\item[{\rm (b)}] For any $\e >0$ there is a trigonometric 
polynomial $b$ on $X$ such that $\|a-b\|_{L^\infty}\leq \e$.
\end{itemize}  
\end{Definition} 

The set $AP(X)$ is a translational invariant unital $C^*$-subalgebra 
of $BC_u(X)$ whose Gelfand spectrum is called {\it the Bohr group} (denoted by $T$). 
All continuous functions on $X$ which are periodic with respect to 
some closed subgroup $\Gamma$ of $X$ with compact quotient $X/\Gamma$
lie in $AP(X)$, but there are many others.         

We can consider the algebra $\A := \langle VO(X)\cdot AP(X)\rangle$
generated by $VO(X)$ and $AP(X)$.
It is obviously a unital $C^*$-subalgebra of $\BC_u(X)$ containing 
$C_0(X)$ and stable by translations. Its Gelfand spectrum is the disjoint union 
$\SA=X\sqcup(R\times T)$, where $R$ is the part at infinity of the Gelfand spectrum 
of $VO(X)$. The relevant quasi-orbits are 
$\big\{\{\var\}\times T\equiv T\big\}_{\var\in R}$. 
This is by no means a general result; it expresses the fact that $VO(X)$
and $AP(X)$ are {\it asymptotically independent}, see \cite{M3} and 
references therein. Actually $AP(X)$ could be replaced by any $C^*$-algebra of 
{\it minimal functions} \cite{M2, M3}. 

Instead of considering arbitrary elements of this algebra $\A$, let us
concentrate on a simple example. 
Assume for simplicity that $V=0$ and that each component of the magnetic field 
is a product of an element of $VO(X)$ and of an element of $AP(X)$, 
{\it i.e.}~$B_{jk} = C_{jk}D_{jk}$, with $C_{jk} \in VO(X)$ and $D_{jk} \in AP(X)$. 
Let us once again invoke the asymptotic values of the matrix valued
function $C:= \{C_{jk}\}_{j,k =1}^N$: $C_\mathbf{x}$ is an asymptotic value  
if and only if there exists a divergent sequence
$\mathbf x=\{x_n\}_{n\in \N}$ such that 
$\sup_{j,k}|(C(x_n)-C_\mathbf x)_{jk}| \to 0$ as $n \to \infty$.
Then one has
\begin{equation*}
\sigma_{\hbox{\rm \tiny ess}}[H_{\!h\!}(B,0)]=\overline{\bigcup_{\mathbf x}
\sigma[H_{\!h\!}(B_\mathbf x,0)]},
\end{equation*}
where $B_\mathbf x$ is the magnetic field whose components are given by
$(B_\mathbf x)_{jk} := (C_\mathbf x)_{jk}D_{jk}\in AP(X)$.
The asymptotic values taken by $C$ at infinity serves as coupling
constants for the magnetic fields of the asymptotic operators, a 
phenomenon already observed in \cite{M1} for Schr\"odinger operators 
without magnetic field. Similarly, if $B_{jk}=C_{jk}+D_{jk}$, with the same 
assumptions on the functions $C_{jk}$ and $D_{jk}$, the asymptotic operators 
are constructed with the almost periodic magnetic fields 
$\{C_\mathbf x+D\}_{\mathbf x\in R}$.

\subsection{Cartesian anisotropy}

In this paragraph we consider another type of spacial anisotropy, 
which is called {\it Cartesian}.
The algebra $C^{\hbox{\rm \tiny cart}}(X)$ consists in the set of 
all continuous functions on $X$ that can be extended to a hypercube 
compactifying $X$. We refer to \cite{R} for a precise definition of 
this algebra and for an extensive study of Schr\"odinger operators
related to this anisotropy (in the absence of magnetic field). 
Let us simply mention that the quasi-orbits are hypercubes
of lower dimensions. We shall restrict here our investigation to a single example 
in the space $\R^2$. In this situation, the set of quasi-orbits consist in 4 
closed segments and 4 points (corners).

For $N=2$, the magnetic field has only one
component $B$ orthogonal to the space $\R^2$. 
Let us assume for simplicity that $B(x_1,x_2)=B_1(x_1)B_2(x_2) + B_0(x_1,x_2)$, where 
$B_0$ belongs to $C_0(\R^2)$ and $B_j(x_j)\to b_j^\pm \in \R$ as 
$x_j \to \pm \infty$. Let also $V$ be of the form $V(x_1,x_2)=V_1(x_1)V_2(x_2) + 
V_0(x_1,x_2)$, where $V_0$ belongs to $C_0(\R^2)$ and $V_j(x_j)\to v_j^\pm \in \R$ as 
$x_j \to \pm \infty$. Then one has
\begin{eqnarray*}
\sigma_{\hbox{\rm \tiny ess}}[H_{\!h\!}(B,V)] & = &
\sigma[H_{\!h\!}(b^-_2B_1,v^-_2V_1)] \;\cup\; \sigma[H_{\!h\!}(b^+_2B_1,v^+_2V_1)] \\
&& \cup\; \sigma[H_{\!h\!}(b^-_1 B_2,v^-_1 V_2)]\;\cup\;\sigma[H_{\!h\!}(b^+_1 B_2,v^+_1 V_2)]. 
\end{eqnarray*}
We stress that each asymptotic operator has a magnetic field that 
depends only on one variable. 
This kind of two dimensional magnetic Schr\"odinger operators was
studied in \cite{I} and \cite{MP4} and exhibits a band spectrum.

\section{Appendix: Some technical results}\label{subsectech} 

Let us now state and prove the auxiliary technical results used in 
the proof of the affiliation criterion.

\begin{Lemma}\label{lemsurgamma}
Assume that the components of the magnetic field $B$ belong to $\A\cap BC^\infty(X)$. 
Then $\gB$ belongs to $\CC^\infty_{\hbox{\rm \tiny pol}}(X \times X; \A)$,
or more precisely:
\begin{enumerate}
\item[{\rm (a)}] for each $x,y \in X$, 
$\;\!\gB(\cdot;x,y) \in \A$,
\item[{\rm (b)}] for each $\alpha, \beta \in \N^N$, there exist
$c>0$, $s_1\geq 0$ and $s_2\geq 0$ 
such that for all $q,x,y \in X$:
\begin{equation*}
\big|\partial_x^\alpha\;\!\partial_y^\beta\;\!
\gB(q;x,y)\big|\leq c \;\!\langle x \rangle^{s_1}
\;\!\langle y \rangle^{s_2}.
\end{equation*}
\end{enumerate}
\end{Lemma}

\begin{proof}
We use the explicit parametrized form of $\gB$
\begin{equation}\label{champm}
\gB(q;x,y)=\exp\Big\{-i\sum_{j,k=1}^Nx_j\;\!y_k
\int_0^1 \de t \int_0^1 \de s\;\!s\;\!B_{jk}\big(q -\hbox{$\frac{1}{2}$}x
-\hbox{$\frac{1}{2}$}y+ sx + st(y-x)\big)\Big\}\;.
\end{equation}
A careful examination of \eqref{champm} leads directly to the results (a) and (b). 
See also the proof of Lemma 4.2 in \cite{MPR}.
\end{proof}

\begin{Lemma}\label{lemsurF}
For each $j \in \{1,\dots,N\}$, each $\alpha, \beta, \gamma \in \N^N$
and each $\mu>\max\{1,s\}$ there exists $c>0$ such that  
\begin{equation}\label{dom}
\big|\partial^\alpha_p \;\! \partial^\beta_k \;\!
\partial^\gamma_l \;\!F_{a,j}(p;k,l)\big|
\leq c \;\! a^{-1/\mu}\;\! \langle p \rangle^{s/\mu-1-|\alpha|}
\;\!\langle k \rangle^s \;\!\langle l \rangle^{2s}
\end{equation}
for all $p,k,l \in X^\star$ and $a \geq -\inf h+1$. 
\end{Lemma}

\begin{proof}
It is enough to show that the expression
\begin{equation}\label{mdg}
\sup_{t \in [0,1]}
\Big|\partial^\alpha_p\;\!\partial^\beta_k \;\!\partial^\gamma_l
\big[(l_j-k_j)\;(\partial_j h)\big(p + (t-1)l-tk\big)\;h_a^{-1}(p-l)
\big]\Big|
\end{equation}
is dominated by the r.h.s.~term of \eqref{dom} with a constant $c$  
not depending on $p$, $k$, $l$ and $a$.

It is easy to see that for any $\delta \in \N^N$, we have 
$\partial^\delta h_a^{-1} = h_a^{-1}\;\! u_{a,\delta}$, where 
$u_{a,\delta} \in S^{-|\delta|}(X^\star)$ uniformly in $a$. By using this, 
the Leibnitz formula and the inequality 
$\langle x+y \rangle \leq \sqrt{2}\langle x\rangle \langle y \rangle$, 
it follows straightforwardly that \eqref{mdg}
is dominated by $c_{\hbox{\tiny 1}}h_a^{-1}(p-l)
\langle p \rangle^{s-1-|\alpha|} \langle k \rangle^{s}
\langle l \rangle^s$ for some $c_{\hbox{\tiny 1}}>0$ independent
of $p$, $k$, $l$ and $a$.
Furthermore, by using the ellipticity of $h$, we see that there exist 
$c_{\hbox{\tiny 2}}>0$ 
and $c_{\hbox{\tiny 3}}>0$ independent of $p,l$ and $a$ such that 
$h_a^{-1}(p-l) \leq c_{\hbox{\tiny 2}}\;\!\langle l \rangle^s
[a+c_{\hbox{\tiny 3}}\;\!\langle p\rangle^s]^{-1}$ for all $p,l \in 
X^\star$.
The final step consists in taking into account the inequality
$a + c_{\hbox{\tiny 3}}\langle p \rangle^{s} \geq 
\mu^{1/\mu}\;\!(\nu c_{\hbox{\tiny 3}})^{1/\nu}\;\!a^{1/\mu}\;\!
\langle p \rangle^{s/\nu}$,
valid for any $\mu \geq 1$, $\nu \geq 1$ with $\mu^{-1}+\nu^{-1}=1$.
\end{proof}

In order to state the next lemma in its full generality, we need the 
definition: 

\begin{Definition}\label{quartiel}
For $s\in \R$, $S^s(X^\star;\A)$ denotes the set of all functions 
$f:X \times X^\star \to \C$ that satisfy: 
\begin{itemize}
\item[{\rm (i)}] $\;\!f(\cdot; p) \in \A$ for all $p \in X^\star$,
\item[{\rm (ii)}] 
$f(q;\cdot)\in C^\infty(X^\star)$, $\forall q\in X$, and for each $\alpha \in \N^N$  
\begin{equation*}
\sup_{q \in X}\|f(q;\cdot )\|_{s,\alpha}:=\sup_{q \in X}\;\! 
\sup_{p \in X^\star}\;\! 
\left[\langle p\rangle^{-s+|\alpha|} \;\!|\partial^\alpha_p f(q;p)|\right]<\infty\;.
\end{equation*}  
\end{itemize}
\end{Definition} 
It is easily seen that the algebraic tensor product $\A \odot S^s(X^\star)$
is contained in $S^s(X^\star; \A)$. 

\begin{Lemma}\label{lemdeRadu}
Let $f$ be an element of $S^t(X^\star;\A)$ with $t<0$. 
Then its partial Fourier transform  $\FF^{-1}(f)$  
is an element of $L^1(X;\A)$ that satisfies for a suitable large integer $m$
\begin{equation}\label{mor}
\|\FF^{-1}(f)\|_{L^1(X;\A)}\leq c\;\!\max_{|\alpha|\leq m}\;\!
\sup_{q \in X}\;\! \|f(q; \cdot)\|_{t,\alpha}\;.
\end{equation} 
\end{Lemma}

\begin{proof}[Proof] 
This is a straightforward adaptation of the proof of \cite[Prop.~1.3.3]{ABG}
(see also \cite[Prop.~1.3.6]{ABG}). We decided to present it in order to put 
into evidence the explicit bound \eqref{mor}. Actually, the arguments needed 
to control the behaviour in the variable $q$ are easy and we leave them to the 
reader; we take simply $f\in S^t(X^\star)$.

Since the case $t\le -N$ is rather simple, we shall concentrate on the more difficult
one: $-N<t<0$. Let us first choose a cutoff function $\chi\in C^\infty_c(X)$ that 
is 1 in a neighbourhood of $0$. One has the estimates:
\begin{eqnarray*}
&& \|(1-\chi)\F^{-1}(f)\|_{L^1} \ \leq \  
C\sum_{|\alpha|=m}\||Q|^{-2m}(1-\chi)\F^{-1}(\partial^{2\alpha} f)\|_{L^1} \\
& \leq & C\Big(\int_X \de x\;\big(1-\chi(x)\big)^2|x|^{-4m}\Big)^{1/2}
\sum_{|\alpha|=m}\|\partial^{2\alpha} f\|_{L^2} \\
& \leq & C'\Big(\int_X \de x\;\big(1-\chi(x)\big)^2|x|^{-4m} \Big)^{1/2}
\ \Big(\int_{X^{\star}}\de p\; \langle p \rangle^{2(t-2m)} \Big)^{1/2} \ 
\underset{|\alpha|=2m}{\max}\|f\|_{t,\alpha}\;,
\end{eqnarray*}
where we take $m \in \N$ with $4m>N$ to make the integrals convergent. 

We study now the behaviour of $\F^{-1}(f)$ near the origin, a more difficult matter. 
Let us fix a second cutoff function $\varphi\in C^\infty(X^\star)$ such that 
$0\leq\varphi\leq 1$, $\varphi(p)=0$ for $|p|\leq 1$ and $\varphi(p)=1$ for $|p|\geq 2$. 
For $b>0$ we set $\varphi_b(p):=\varphi(bp)$. We have:
\begin{equation*}
\big|\big\{\F^{-1}\big((1-\varphi_b)f\big)\big\}(y)\big|\;\leq
\int\limits_{|p|\leq2/b}\de p\;|f(p)| \; \leq \; \|f\|_{t,0}
\int\limits_{|p|<2/b}\de p\;|p|^{t} \leq C\;\|f\|_{t,0}\;b^{-N-t}.
\end{equation*}
Moreover, if $m\in 2\N$ with $m\geq N+1$, then one has:
\begin{eqnarray*}
&& |y|^{m}|[\F^{-1}(\varphi_b f)](y)|\ \leq \
C\sum_{|\alpha|=m}\big|\big[\F^{-1}\big(\partial^\alpha(\varphi_b f)\big)\big](y)\big| \\
& \leq & C\sum_{|\alpha|=m}\;\sum_{\beta\leq\alpha}C_\alpha^\beta \;b^{|\alpha-\beta|}
\int_{X^\star}\de p\;|(\partial^{\alpha-\beta}\varphi)(bp)|\;|(\partial^\beta f)(p)| \\
& \leq & C'\underset{|\alpha|\leq m}{\max}\|f\|_{t,\alpha}\;
\Big\{\int\limits_{|p|\geq 1/b}\de p\;|p|^{t-m} + 
\sum_{|\beta|< m}b^{m-|\beta|}\int\limits_{1/b<|p|<2/b}\de p\;|p|^{t-|\beta|}\Big\} \\
& = & C''\;\underset{|\alpha|\leq m}{\max}\|f\|_{t,\alpha}\;b^{m-N-t}.
\end{eqnarray*}
By fixing $b:=|y|$, we get
$|[\F^{-1}(\varphi_{|y|} f)](y)|\leq C''\;\underset{|\alpha|\leq m}{\max}\|f\|_{t,\alpha}
\;|y|^{-N-t}$.
The singularity at the origin is integrable, and putting all the inequalities together 
we obtain \eqref{mor}.
\end{proof}

\section*{Acknowledgment}
Serge Richard is supported by the Swiss National Science Foundation. Marius Mantoiu and Radu Purice acknowledge partial support from the CERES Contract No. 3-28/2003 with the Romanian Ministery of Education and Research. Part of this work has been completed while the authors visited the University of Geneva and we express our gratitude to Prof. Werner Amrein for his kind hospitality and the stimulating discusions.

\end{document}